\documentclass[english,envcountsect,envcountsame]{svjour3}
\usepackage[T1]{fontenc}
\usepackage[latin9]{inputenc}
\usepackage{babel}
\usepackage{refstyle}
\usepackage{url}
\usepackage{amsmath}
\usepackage{amssymb}
\usepackage{graphicx}
\usepackage[numbers,sort&compress]{natbib}
\usepackage[unicode=true,pdfusetitle,
 bookmarks=true,bookmarksnumbered=false,bookmarksopen=false,
 breaklinks=false,pdfborder={0 0 0},pdfborderstyle={},backref=false,colorlinks=false]
 {hyperref}

\makeatletter

\AtBeginDocument{\providecommand\exaref[1]{\ref{exa:#1}}}
\AtBeginDocument{\providecommand\figref[1]{\ref{fig:#1}}}
\AtBeginDocument{\providecommand\thmref[1]{\ref{thm:#1}}}
\AtBeginDocument{\providecommand\corref[1]{\ref{cor:#1}}}
\AtBeginDocument{\providecommand\propref[1]{\ref{prop:#1}}}
\AtBeginDocument{\providecommand\lemref[1]{\ref{lem:#1}}}
\AtBeginDocument{\providecommand\secref[1]{\ref{sec:#1}}}
\AtBeginDocument{\providecommand\Lemref[1]{\ref{Lem:#1}}}
\AtBeginDocument{\providecommand\appref[1]{\ref{app:#1}}}
\RS@ifundefined{subsecref}
  {\newref{subsec}{name = \RSsectxt}}
  {}
\RS@ifundefined{thmref}
  {\def\RSthmtxt{theorem~}\newref{thm}{name = \RSthmtxt}}
  {}
\RS@ifundefined{lemref}
  {\def\RSlemtxt{lemma~}\newref{lem}{name = \RSlemtxt}}
  {}

\numberwithin{equation}{section}

\usepackage{xcolor}

\usepackage{url}

\smartqed
\newref{prop}{name=Proposition~,Name=Proposition~}
\newref{thm}{name=Theorem~,Name=Theorem~}
\newref{lem}{name=Lemma~,Name=Lemma~}
\newref{cor}{name=Corollary~,Name=Corollary~}
\newref{sec}{name=Section~,Name=Section~}
\newref{app}{name=Appendix~,Name=Appendix~}
\newref{exa}{name=Example~,Name=Example~}
\newref{fig}{name=Fig.~,Name=Fig.~}

\makeatother

\begin{document}
\title{Improved Global Guarantees for the Nonconvex Burer--Monteiro Factorization
via Rank Overparameterization\thanks{Financial support for this work was provided in part by the NSF CAREER
Award ECCS-2047462 and in part by C3.ai Inc. and the Microsoft Corporation
via the C3.ai Digital Transformation Institute.}}
\titlerunning{Improved Global Guarantees via Rank Overparameterization}
\author{Richard Y. Zhang\\
}
\institute{Dept. of Electrical and Computer Engineering\\
University of Illinois at Urbana-Champaign\\
306 N Wright St, Urbana, IL 61801\\
\email{ryz@illinois.edu}}
\date{~}
\maketitle
\begin{abstract}
We consider minimizing a twice-differentiable, $L$-smooth, and $\mu$-strongly
convex objective $\phi$ over an $n\times n$ positive semidefinite
matrix $M\succeq0$, under the assumption that the minimizer $M^{\star}$
has low rank $r^{\star}\ll n$. Following the Burer--Monteiro approach,
we instead minimize the nonconvex objective $f(X)=\phi(XX^{T})$ over
a factor matrix $X$ of size $n\times r$. This substantially reduces
the number of variables from $O(n^{2})$ to as few as $O(n)$ and
also enforces positive semidefiniteness for free, but at the cost
of giving up the convexity of the original problem. In this paper,
we prove that if the search rank $r\ge r^{\star}$ is overparameterized
by a \emph{constant factor} with respect to the true rank $r^{\star}$,
namely as in $r>\frac{1}{4}(L/\mu-1)^{2}r^{\star}$, then despite
nonconvexity, local optimization is guaranteed to globally converge
from any initial point to the global optimum. This significantly improves
upon a previous rank overparameterization threshold of $r\ge n$,
which we show is sharp in the absence of smoothness and strong convexity,
but would increase the number of variables back up to $O(n^{2})$.
Conversely, without rank overparameterization, we prove that such
a global guarantee is possible if and only if $\phi$ is almost perfectly
conditioned, with a condition number of $L/\mu<3$. Therefore, we
conclude that a small amount of overparameterization can lead to large
improvements in theoretical guarantees for the nonconvex Burer--Monteiro
factorization.
\end{abstract}

\section{Introduction}

\global\long\def\R{\mathbb{R}}%
\global\long\def\S{\mathbb{S}}%
\global\long\def\AA{\mathcal{A}}%
\global\long\def\HH{\mathcal{H}}%
\global\long\def\II{\mathcal{I}}%
\global\long\def\M{\mathcal{M}}%
\global\long\def\rank{\mathrm{rank}}%
\global\long\def\av{\mathrm{av}}%
\global\long\def\poly{\mathrm{poly}}%
\global\long\def\nnz{\mathrm{nnz}}%
\global\long\def\e{\mathbf{e}}%
\global\long\def\J{\mathbf{J}}%
\global\long\def\A{\mathbf{A}}%
\global\long\def\H{\mathbf{H}}%
\global\long\def\G{\mathbf{G}}%
\global\long\def\vect{\mathrm{vec}}%
\global\long\def\mat{\mathrm{mat}}%
\global\long\def\ub{\mathrm{ub}}%
\global\long\def\lb{\mathrm{lb}}%
\global\long\def\tr{\mathrm{tr}}%
\global\long\def\eqdef{\overset{\mathrm{def}}{=}}%
\global\long\def\LMI{\mathrm{LMI}}%
\global\long\def\d{\mathrm{d}}%
\global\long\def\one{\mathbf{1}}%
\global\long\def\diag{\mathrm{diag}}%
\global\long\def\spur{\mathrm{spur}}%
\global\long\def\half{\frac{1}{2}}%
\global\long\def\f{\mathbf{f}}%
\global\long\def\Q{\mathbf{Q}}%
\global\long\def\inner#1#2{\left\langle #1,#2\right\rangle }%
Consider minimizing a convex objective $\phi$ over an $n\times n$
positive semidefinite matrix $M\succeq0$, as in
\begin{equation}
M^{\star}=\quad\mathrm{minimize}\quad\phi(M)\quad\text{over }\quad M\succeq0,\tag{P}\label{eq:cvx}
\end{equation}
In practice, the matrix order $n$ is often so large---from tens
of thousands to hundreds of millions---that even explicitly forming
the $n^{2}$ elements of the $n\times n$ matrix variable $M$ would
be intractable. Where $r^{\star}=\rank(M^{\star})$ is known a priori
to be small, the standard approach is the nonconvex Burer--Monteiro
factorization~\citep{burer2003nonlinear}, which rewrites $M=XX^{T}$
and then directly optimizes over $X$ using a local optimization algorithm:
\begin{equation}
X^{\star}=\quad\mathrm{minimize}\quad f(X)\eqdef\phi(XX^{T})\text{ where }X\text{ is }n\times r\text{ and }r\ge r^{\star}.\tag{BM}\label{eq:ncvx}
\end{equation}
This way, the semidefinite constraint $M=XX^{T}\succeq0$ is automatically
enforced for free. Moreover, by using a small search rank $r\ll n$,
the number of variables is reduced to $O(n)$. However, this improved
scalability comes at the cost of giving up the convexity of (\ref{eq:cvx}).
In principle, local optimization can fail by getting stuck at a \emph{spurious}
local minimum of (\ref{eq:ncvx})---a local minimum that is strictly
worse than that of the global minimum. In practice, however, there
seems to be many situations where this failure mode does not occur~\citep{candes2015phase,rosen2019se,zhang2019spurious}.

In pursuit of a theoretical explanation for the empirical effectiveness
of the Burer--Monteiro approach, it was discovered that if the convex
function $\phi$ is sufficiently \emph{well-conditioned}, then $f$
has \emph{no spurious local minima}~\citep{bhojanapalli2016global,ge2017nospurious,park2018finding,zhu2018global,li2019non}.
To state this precisely, suppose that $\phi$ is twice-differentiable,
$L$-smooth, and $\mu$-strongly convex over the entire set of $n\times n$
real symmetric matrices $\S^{n}$, as in
\[
\mu\|E\|_{F}^{2}\le\inner{\nabla^{2}\phi(M)[E]}E\le L\|E\|_{F}^{2}\qquad\text{for all }M\in\S^{n},
\]
where $\inner EF\eqdef\tr(E^{T}F)$ and $\|E\|_{F}\eqdef\sqrt{\inner EE}$
denote the matrix Euclidean inner product and norm respectively. \citet{bhojanapalli2016global}
pointed out that if $\phi$ has condition number $L/\mu<3/2$, then
the usual second-order necessary conditions for \emph{local} optimality
in $f$ are also sufficient for \emph{global} optimality, meaning
that
\[
\nabla f(X)=0,\quad\nabla^{2}f(X)\succeq0\quad\iff\quad f(X)=\min_{U}f(U).
\]
Moreover, satisfying these conditions to $\epsilon$-accuracy will
yield a point within a $\rho$-neighborhood of a globally optimal
solution. Where this guarantee holds, local optimization cannot fail
by getting stuck, because every local minimum is also a global minimum,
and every saddle-point has sufficiently negative curvature to allow
escape. (See also \citep{ge2017nospurious,park2018finding,zhu2018global,li2019non}
for improvements and extensions on this result.)

While powerful, this global guarantee is unfortunately also very conservative.
Practical choices of $\phi$ can be assumed to have a \emph{finite}
condition number $L/\mu$, because they tend to be $L$-smooth by
formulation, and $\mu$-strong convexity can be made to hold by adding
a small regularizer, as in $\phi_{\mu}(M)=\phi(M)+\frac{\mu}{2}\|M\|_{F}^{2}$.
But adding a regularizer would also worsen the accuracy and quality
of the solution, so it is unrealistic to assume condition numbers
as small as $L/\mu<3/2$. On the other hand, \citet{zhang2018much}
stated a counterexample $\phi$ with $L/\mu=3$ but whose $f$ admits
a spurious second-order point (see also~\citep{zhu2018global,zhang2019sharp,li2019non}).
For more realistic condition numbers $L/\mu\ge3$, the Burer--Monteiro
approach may continue to work well in practice, but one can no longer
rule out failure by getting stuck at a spurious local minimum.

\begin{figure}
\includegraphics[width=0.5\columnwidth]{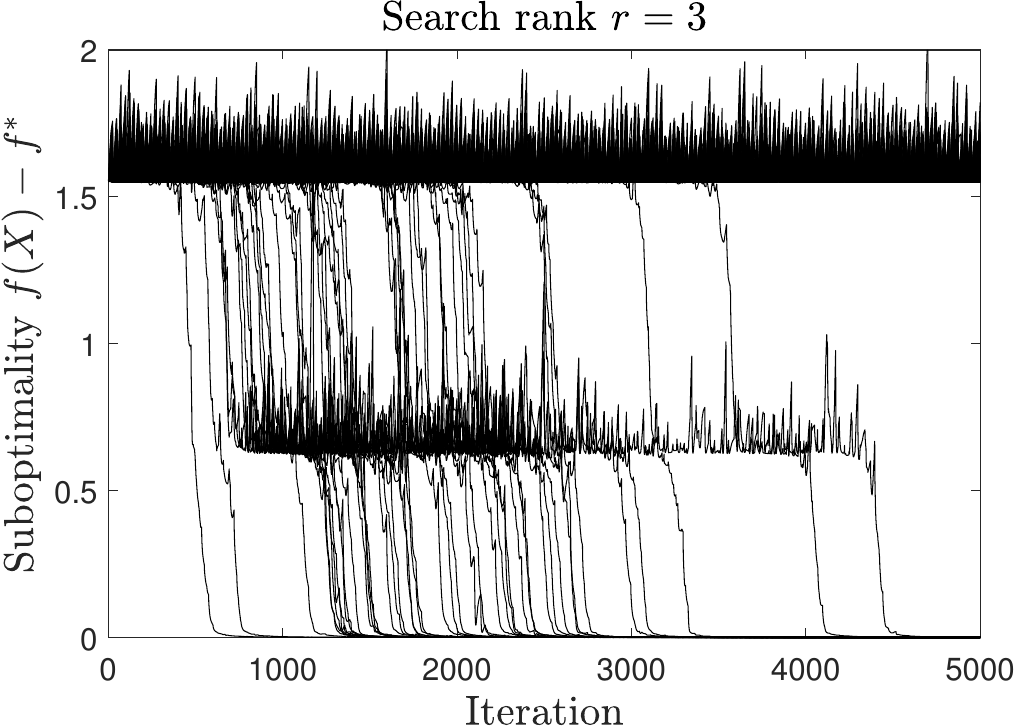}\hspace*{\fill}\includegraphics[width=0.5\columnwidth]{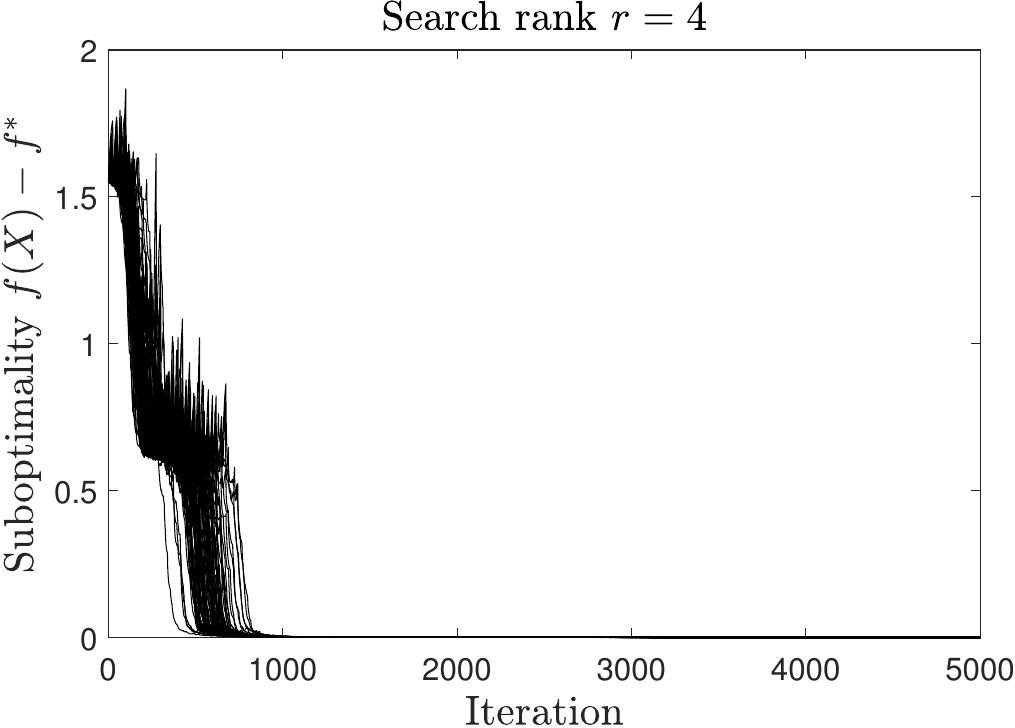}

\caption{\label{fig:overparam}\textbf{Overparameterization eliminates spurious
local minima.} Stochastic gradient descent (SGD) with Nesterov momentum~\citep{sutskever2013importance}
applied to an $f(X)\protect\eqdef\phi(XX^{T})$ with a spurious second-order
point $X_{\mathrm{spur}}$ for $r=3$: \textbf{(Left)} With search
rank $r=3$, GD remains stuck at $X\approx X_{\protect\spur}$, resulting
in 55 failures out of 100 trials. \textbf{(Right)} Overparameterizing
to $r=4$ eliminates $X_{\protect\spur}$ as a spurious second-order
point, and GD now succeeds in all 100 trials. (Set $\phi(M)=\sum_{i,j=1}^{n}\phi_{i,j}(M)$
where $\phi_{i,j}(M)=\frac{1}{2}|\protect\inner{A^{(i,j)}}{M-M^{\star}}|^{2}$
as in \exaref{overparam} with $n=5$, $r=3,$ and $r^{\star}=2$,
set $V=0$ and uniformly sample $X$ from $\|X-X_{\mathrm{spur}}\|_{F}\le0.1$,
and then run $V_{\mathrm{new}}=\beta V-\alpha\nabla f_{i,j}(X)$ and
$X_{\mathrm{new}}=X+\beta V_{\mathrm{new}}-\alpha\nabla f_{i,j}(X)$,
with learning rate $\alpha=1\times10^{-1}$ and momentum $\beta=0.9$.
Sample indices $i,j$ are randomly shuffled every 1 epoch = 25 iterations.)}
\end{figure}

\subsection{Main result}

Surprisingly, we show in this paper that \emph{overparameterizing}
the search rank $r>r^{\star}$ can improve existing global guarantees
to all finite values of the condition number $L/\mu$---far beyond
the apparently fundamental barrier of $L/\mu<3$. We are inspired
by parallel work on the Burer--Monteiro approach for semidefinite
programming (see our literature review in Section~\ref{subsec:Overparameterization}),
where it was discovered that progressively overparameterizing the
search rank $r$ makes the nonconvexity increasingly benign. The same
empirical observation is easily replicated in our unconstrained setting:
for a fixed unconstrained convex $\phi$, the corresponding nonconvex
$f$ admits progressively fewer spurious local minima as its search
rank $r$ is progressively increased past $r^{\star}$ (see \figref{overparam}).
Once the overparameterization ratio $r/r^{\star}$ exceeds some constant
threshold, local optimization consistently succeeds at globally minimizing
$f$. 

Our main result rigorously justifies the above empirical observation:
so long as the convex function $\phi$ has a bounded condition number
$L/\mu=O(1)$ with respect to $n$, the corresponding nonconvex function
$f$ is guaranteed to have no spurious local minima for a sufficiently
large search rank $r=O(r^{\star})$. To the best of our knowledge,
this is the first rigorous proof that a \emph{constant-factor} overparameterization
ratio $r/r^{\star}$ can eliminate spurious local minima.
\begin{theorem}[Overparameterization]
\label{thm:main}Let $\phi:\S^{n}\to\R$ be twice-differentiable,
$L$-smooth and $\mu$-strongly convex, let the minimizer $M^{\star}=\arg\min_{M\succeq0}\phi(M)$
have true rank $r^{\star}=\rank(M^{\star})$.
\begin{itemize}
\item (Sufficiency) If $r>\frac{1}{4}(L/\mu-1)^{2}\,r^{\star}$ and $r^{\star}\le r<n$,
then the Burer--Monteiro function $f:\R^{n\times r}\to\R$ defined
as $f(U)\eqdef\phi(UU^{T})$ has no spurious local minima:
\[
\nabla f(X)=0,\quad\nabla^{2}f(X)\succeq0\quad\iff\quad f(X)=\min_{U}f(U).
\]
\item (Necessity) If $r\le\frac{1}{4}(L/\mu-1)^{2}\,r^{\star}$ and $r^{\star}\le r<n$,
then there exists an $L$-smooth and $\mu$-strongly convex quadratic
counterexample $\phi_{0}$ whose minimizer $M_{0}^{\star}=\arg\min_{M\succeq0}\phi_{0}(M)$
has $\rank(M_{0}^{\star})=r^{\star}$, but whose $f_{0}(U)\eqdef\phi_{0}(UU^{T})$
admits an $n\times r$ spurious second-order point $X_{\spur}$:
\[
\nabla f_{0}(X_{\spur})=0,\quad\nabla^{2}f_{0}(X_{\spur})\succeq0,\quad f_{0}(X_{\spur})-\min_{U}f_{0}(U)>\frac{\mu}{2}\cdot\|M_{0}^{\star}\|_{F}^{2}.
\]
\end{itemize}
\end{theorem}
Once the search rank is large enough to satisfy $r\ge n$, one can
easily adapt existing results to prove that $f$ has no spurious local
minima; see \citet[Corollary 8]{journee2010low} and \citet[Corollary 3.2]{boumal2020deterministic}.
Note that this guarantee does not require $\phi$ to be strongly convex,
nor uniformly $L$-smooth over its entire domain. Of course, setting
$r\ge n$ would also force us to optimize over $O(n^{2})$ matrix
elements in $X$, thereby obviating the computational advantages of
the Burer--Monteiro factorization in the first place. 

In fact, the necessary condition in \thmref{main} shows that the
overparameterization rank threshold $r\ge n$ is sharp. Without smoothness
and strong convexity, it is inherently impossible to make a global
guarantee for a search rank $r<n$, due to existence of a counterexample
(\exaref{overparam}). The sufficient condition in \thmref{main}
is able to improve upon the sharp threshold of $r\ge n$ only because
it further keeps the condition number $L/\mu$ bounded. Our key insight
is that highly overparameterized counterexamples $\phi$ do exist,
but their condition numbers $L/\mu$ must necessarily diverge to infinity.
If instead the condition number $L/\mu$ can be kept constant---which
is indeed the case if the convex function $\phi$ is held fixed---then
we can progressively increase the overparameterization ratio $r/r^{\star}$,
and use the inexistence of a well-conditioned overparameterized counterexample
to prove that $f$ has no spurious local minima. 

Without explicitly accounting for overparameterization, the best global
guarantee that we can make is for almost-perfect choices of $\phi$
with conditions numbers $L/\mu<3$. Indeed, we obtain the following
by allowing the true rank $r^{\star}$ in \thmref{main} to take on
all values smaller than the search rank $r$ (i.e. by substituting
the upper-bound $r^{\star}\le r$). It provides an affirmative answer
to a conjecture first posed in \citet{zhang2018much} that is widely
believed to be true in the existing literature~\citep{zhu2018global,zhang2019sharp,li2019non}.
\begin{corollary}[Exact parameterization]
\label{cor:half}Let $\phi:\S^{n}\to\R$ be twice-differentiable,
$L$-smooth and $\mu$-strongly convex, and let the search rank $r$
satisfy $r\ge\rank(M^{\star})$ where $M^{\star}=\arg\min_{M\succeq0}\phi(M)$. 
\begin{itemize}
\item If $L/\mu<3$, then $f(U)\eqdef\phi(UU^{T})$ for $U\in\R^{n\times r}$
has no spurious local minima. 
\item If $L/\mu\ge3$, then nothing can be said due to the existence of
a counterexample.
\end{itemize}
\end{corollary}
Comparing \thmref{main} and \corref{half}, we conclude that a small
amount of overparameterization can lead to large improvements in global
guarantees based on the condition number $L/\mu$. As an important
future work, we expect that explicitly accounting for rank overparameterization
should also greatly improve other theoretical guarantees, including
local guarantees based on a specific initial point and choice of algorithm~\citep{zhang2019sharp,zhang2020many},
and also application-specific global guarantees, based on narrower
structures like incoherence~\citep{ge2016matrix,ge2017nospurious,ma2019implicit},
low-rank measurements~\citep{chen2015solving,candes2015phase,chen2019gradient},
and structured sparsity~\citep{molybog2020role}. In this regard,
an important contribution of this paper is to provide a set of tools
to reason about the inexistence of overparameterized counterexamples.

Finally, the global guarantee in \thmref{main}---though improved---is
unfortunately still quite conservative. For example, with a condition
number of $L/\mu\approx200$, the overparameterization ratio $r/r^{\star}=\frac{1}{4}(L/\mu-1)^{2}\approx10^{4}$
required by \thmref{main} would be far too large for practical use.
On the other hand, this conservatism appears fundamental; \thmref{main}
is already sharp, and as we explain in Section~\ref{subsec:Matrix-sensing}
below, cannot be improved by adopting the rank-restricted versions
of $\mu$-strong convexity and $L$-smoothness. As explained earlier,
it is generally undesirable to use a regularizer to improve the condition
number, as it would also worsen accuracy. A better approach is to
adopt a preconditioner, as is typically done to improve the convergence
rate of algorithms, but this is not a general-purpose approach and
would require further investigations. It also remains an important
open question whether there exists even stronger structures than $L$-smoothness
and $\mu$-strong convexity, that can better control the benign nonconvexity
of the Burer--Monteiro factorization. 

\subsection{Algorithmic implications}

While \thmref{main} guarantees global optimality for a point $X$
that \emph{exactly} satisfies the second-order conditions for local
optimality, as in $\nabla f(X)=0$ and $\nabla f(X)\succeq0$, practical
algorithms like gradient descent~\citep{jin2017escape} and trust-region
methods~\citep{nesterov2006cubic,cartis2011adaptive1} are only able
to compute a point $\tilde{X}$ that \emph{approximately} satisfies
these conditions up to some prescribed accuracy $\epsilon>0$, as
in $\|\nabla f(\tilde{X})\|_{F}\le\epsilon$ and $\nabla^{2}f(\tilde{X})\succeq-\epsilon I$.
In order to guarantee global convergence towards global optimality,
we additionally require the \emph{strict saddle property}, which says
that every $\epsilon$ approximate second-order point $\tilde{X}$
is also guaranteed to be $\rho$-close to an exact second-order point
$X$~\citep{ge2015escaping,jin2017escape,ge2017nospurious}. Fortunately
within our setting, this strengthened property comes essentially for
free.
\begin{proposition}[Strict saddle property]
\label{prop:strictsaddle}Let $\phi:\S^{n}\to\R$ be twice-differentiable,
$L$-smooth and $\mu$-strongly convex, and let the search rank $r$
satisfy $r\ge\rank(M^{\star})$ where $M^{\star}=\arg\min_{M\succeq0}\phi(M)$.
If the Burer--Monteiro function $f(U)\eqdef\phi(UU^{T})$ has no
spurious local minima
\[
\nabla f(X)=0,\quad\nabla^{2}f(X)\succeq0\quad\iff\quad f(X)=\min_{U}f(U)
\]
then $f$ also satisfies the strict saddle property: 
\[
\|\nabla f(\tilde{X})\|_{F}\le\epsilon(\delta),\quad\nabla^{2}f(\tilde{X})\succeq-\epsilon(\delta)I\quad\implies\quad f(X)-\min_{U}f(U)<\delta
\]
where $\epsilon$ is a nondecreasing function that satisfies $0<\epsilon(\delta)\le1$
for all $\delta>0$.
\end{proposition}
Note that \propref{strictsaddle} promises only global convergence
towards global optimality in \emph{finite} time, and not necessarily
in \emph{reasonable} time.  Nevertheless, once convergence can be
rigorously established, we expect that a convergence rate should readily
follow with a more careful analysis. In the existing literature, once
it became clear that $f$ has no spurious local minima with $L/\mu<3/2$,
subsequent refinements quickly showed that  every $\epsilon$ approximate
second-order point $\tilde{X}$ is guaranteed to be $O(\epsilon)$-close
to an exact second-order point, and therefore at most $O(\epsilon^{2})$-globally
suboptimal~\citep{ge2017nospurious,zhu2018global,li2019non}. With
exact rank parameterization $r=r^{\star}$, gradient descent globally
converges at a linear rate to a $\delta$-globally suboptimal point
$\tilde{X}$ satisfying $f(\tilde{X})-\min_{U}f(U)<\delta$ in $O(\log(1/\delta))$
iterations, as if $f$ were strongly convex~\citep{jin2017escape,ge2017nospurious}.
With rank over-parameterization $r>r^{\star}$, gradient descent slows
down to a sublinear rate; instead, \citet{zhang2021preconditioned,zhang2022preconditioned}
recently proposed a preconditioned gradient descent algorithm and
proved that it restores global linear convergence. We expect that
all of these results can be generalized to the overparameterized regime
of $r/r^{\star}>\frac{1}{4}(L/\mu-1)^{2}$, by suitably extending
the proof of \thmref{main} to approximate second-order points. We
highlight this as another important future direction. 

\subsection{Related work}

\subsubsection{\label{subsec:Overparameterization}Rank overparameterization for
semidefinite programs}

Rank overparameterization for the Burer--Monteiro factorization is
currently best understood in the context of solving semidefinite programs
\begin{equation}
\text{minimize}\quad\inner C{XX^{T}}\quad\text{ subject to }\quad\AA(XX^{T})=b,\quad X\text{ is }n\times r,\tag{SDP}\label{eq:sdp}
\end{equation}
in which $\AA:\S^{n}\to\R^{m}$ models a set of $m<\frac{1}{2}n(n+1)$
linear constraints with right-hand side $b\in\R^{m}$, and $C\in\S^{n}$
models a linear cost. The critical insight (which also applies to
our setting) is that a spurious local minimum at search rank $r$
must correspond to a strict saddle-point at search rank $r+1$ that
can be escaped~\citep{journee2010low,boumal2016non,boumal2020deterministic}.
The algorithm obtained by incrementally overparameterizing the search
rank $r$ this way---known as the \emph{staircase} method~\citep{boumal2016non,boumal2020deterministic}---often
terminates at $r\approx10$ with the global solution~\citep{rosen2019se,chiu2022overcoming}.
To provide a correctness guarantee, \citet[Corollary 8]{journee2010low}
proved that (\ref{eq:sdp}) has no spurious local minima if $r\ge n$,
so the staircase method must terminate. \citet{boumal2016non,boumal2020deterministic}
improved this threshold to $r>\sqrt{2m+\frac{1}{4}}-\frac{1}{2}$
for a \emph{generic} instance of (\ref{eq:sdp}), which \citet{waldspurger2020rank}
showed is essentially sharp due to the existence of counterexamples.
More recently, \citet{ocarroll2022burer} established, in the non-generic
worst case, that the threshold $r\ge n$ is essentially sharp.

Although (\ref{eq:sdp}) and (\ref{eq:ncvx}) consider two very different
settings, a connection between them can be established through the
augmented Lagrangian
\[
f_{\beta}(X)\eqdef\phi_{\beta}(XX^{T};y)\quad\text{ where }\phi_{\beta}(M;y)\eqdef\inner{C-\AA^{T}(y)}M+\frac{\beta}{2}\|\AA(M)-b\|^{2}.
\]
In one direction, the original papers of \citet{burer2003nonlinear,burer2005local}
proposed to solve (\ref{eq:sdp}) via this particular instance of
(\ref{eq:ncvx}). However, $\phi_{\beta}$ can never be strongly convex,
so the best global guarantee for $f_{\beta}$ is the same rank threshold
$r\ge n$, which is also sharp within our context via \thmref{main}.
Imposing a strongly convex regularizer, as in $f_{\beta,\mu}(X)\eqdef f_{\beta}(X)+\frac{\mu}{2}\|X^{T}X\|_{F}^{2}$,
could potentially allow global guarantees for $r\ll n$ via \thmref{main}
to be eventually extended to (\ref{eq:sdp}), but this could again
come at the cost of introducing substantial errors. 

In reverse, the counterexamples originally formulated for (\ref{eq:sdp})
can be adapted to (\ref{eq:ncvx}) under second-order sufficiency
conditions. The following lemma establishes this classical connection
(see e.g.~\citep[Section 2.2.2]{bertsekas2014constrained}) while
quotienting out the rotational invariance $f_{\beta}(X)=f_{\beta}(XR)$
for $RR^{T}=I_{r}$.
\begin{lemma}
\label{lem:aug}Let $\mathcal{V}=\{V\in\R^{n\times r}:XV^{T}+VX^{T}\ne0\}$.
If $X\in\R^{n\times r},y\in\R^{m}$ satisfy the second-order sufficient
conditions for local optimality in (\ref{eq:sdp}):\begin{subequations}\label{eq:sdp_soc}
\begin{gather}
[C-\AA^{T}(y)]X=0,\qquad\AA(XX^{T})=b,\label{eq:sdp_soc_1}\\
V\in\mathcal{V},\quad\AA(XV^{T}+VX^{T})=0\quad\implies\quad\inner{C-\AA^{T}(y)}{VV^{T}}>0,\label{eq:sdp_soc_2}
\end{gather}
\end{subequations}then there exists $\beta\ge0$ such that $X$ satisfies
$\nabla f_{\beta}(X)=0,$ $\nabla^{2}f_{\beta}(X)\succeq0,$ and $\inner{\nabla^{2}f_{\beta}(X)[V]}V>0$
for all $V\in\mathcal{V}$. 
\end{lemma}
\begin{proof}
Explicitly write $\nabla f_{\beta}(X)=0$ and $\inner{\nabla^{2}f_{\beta}(X)[V]}V>0$
for $V\in\mathcal{V}$ as\begin{subequations}
\begin{gather}
[C-\AA^{T}(y_{+})]X=0\quad\text{ where }y_{+}=y+\beta[b-\AA(XX^{T})],\label{eq:aug_soc_1}\\
\inner{C-\AA^{T}(y_{+})}{VV^{T}}+\frac{\beta}{2}\|\AA(XV^{T}+VX^{T})\|^{2}>0\text{ for all }V\in\mathcal{V}.\label{eq:aug_soc_2}
\end{gather}
\end{subequations}Clearly, (\ref{eq:sdp_soc_1}) implies (\ref{eq:aug_soc_1})
because $y_{+}=y$. Observe that $\mathcal{V}$ is a linear subspace;
we evoke the homogeneous nonstrict version of the S-lemma~\citep[Proposition~3.2]{polik2007survey}
to conclude that there exists no $V\in\mathcal{V}$ that jointly satisfies
$\|\AA(XV^{T}+VX^{T})\|^{2}\le0$ and $\inner{C-\AA^{T}(y)}{VV^{T}}\le0$
if and only if (\ref{eq:aug_soc_2}) holds for some $\beta\ge0$.
Finally, any $U\notin\mathcal{V}$ can be written as $U=XA$ for $A=-A^{T}$,
and in this case, $\nabla^{2}f_{\beta}(X)[U]=[C-\AA^{T}(y_{+})]XA=0$
via (\ref{eq:aug_soc_1}). Therefore, we conclude that $\nabla^{2}f_{\beta}(X)\succeq0$
holds. \qed
\end{proof}

The second-order sufficient conditions are satisfied by the counterexamples
of \citet{waldspurger2020rank} for the MAXCUT problem (where it is
called \emph{nondegeneracy}; see~\citep[Definition~2.7]{waldspurger2020rank}),
so applying \lemref{aug} results in a counterexample $\phi$ with
true rank $r^{\star}=1$ but a spurious second-order point at $r\ge\sqrt{2n}$.
The counterexamples of \citet{ocarroll2022burer} do not satisfy these
conditions, but nevertheless hint at the existence of a counterexample
$\phi$ with $r^{\star}=1$ and $r=\Omega(n)$. \citet[Theorem~5]{bhojanapalli2018smoothed}
stated a counterexample $\phi$ with $r^{\star}=1$ and $r=n-1$ that
fills this gap. Finally, \thmref{main} states the optimal counterexample
$\phi^{\star}$ with $r^{\star}=1$ and $r=n-1$, and whose finite
condition number $L/\mu=2\sqrt{n-1}+1$ is in fact as small as possible.

\subsubsection{\label{subsec:Matrix-sensing}Matrix sensing}

Our results are closely related to a problem known in the literature
as \emph{matrix sensing}, which  considers a nonlinear least-squares
objective like $f(X)=\frac{1}{2}\|\AA(XX^{T}-M^{\star})\|^{2}$, in
which  the measurement operator $\AA$ is assumed to satisfy the
$(\delta,r+r^{\star})$-\emph{restricted isometry property} or $(\delta,r+r^{\star})$-RIP,
with constant $\delta\in[0,1)$:
\[
(1-\delta)\|E\|_{F}^{2}\le\|\AA(E)\|^{2}\le(1+\delta)\|E\|_{F}^{2}\qquad\text{ for all }\rank(E)\le r+r^{\star}.
\]
The matrix sensing problem has drawn considerable interest because
it gives some of the most convincing and unambiguous demonstration
of a nonconvex problem with a ``benign landscape''~\citep{bhojanapalli2016global,ge2017nospurious,chi2019nonconvex}.
In practice, applying simple gradient descent to $f$ consistently
results in rapid linear convergence to the global optimum, as if it
were strongly convex~\citep{zheng2015convergent,tu2016low,jin2017escape}.
Rigorously, it is known that $\delta<1/5$ is sufficient~\citep{bhojanapalli2016global,ge2017nospurious,park2018finding,li2019non,zhu2018global}
and $\delta<1/2$ is necessary~\citep{zhang2018much,zhu2018global,li2019non}
for the nonconvex function $f$ to satisfy the strict saddle property
and have no spurious local minima. 

A twice-differentiable, $L$-smooth, and $\mu$-strongly convex function
$\phi$ by definition satisfies $(\delta,r+r^{\star})$-RIP with $\delta=\frac{L-\mu}{L+\mu}$,
so existing results can be adapted to imply that a condition number
of $L/\mu<3/2$ is sufficient and $L/\mu<3$ is necessary for spurious
local minima not to exist. However, note that RIP is a more general
assumption than smoothness and strong convexity, so our global guarantee
in \thmref{main} does not imply a similar result under RIP. Fortunately,
its proof can be strengthened to cover this more general case, although
it does not come with any substantial improvements.
\begin{corollary}[Restricted isometry property]
\label{cor:RIP}Let $M^{\star}\succeq0$ have true rank $r^{\star}=\rank(M^{\star})$
and let $\AA:\S^{n}\to\R^{m}$ satisfy $(\delta,r+r^{\star})$-RIP.
\begin{itemize}
\item (Sufficiency) If $\delta<1/(1+\sqrt{r^{\star}/r})$, then $f(X)\eqdef\|\AA(XX^{T}-M^{\star})\|^{2}$
has no spurious local minima. 
\item (Necessity) If $\delta\ge1/(1+\sqrt{r^{\star}/r})$, then there exists
a counterexample $\AA_{0}$ that satisfies $(\delta,r+r^{\star})$-RIP,
but whose $f_{0}(X)\eqdef\|\AA_{0}(XX^{T}-M^{\star})\|^{2}$ admits
a spurious second-order point.
\end{itemize}
\end{corollary}
\begin{proof}
Repeat the proof of \thmref{main} to minimize the function $\delta(X,Z)$
as posed in \citep[Theorem 8]{zhang2019sharp}. We defer the specific
details of the sufficient condition to the technical report~\citep{zhang2021sharp},
and prove necessity using \exaref{overparam} in \secref{counterexamples}.
\qed
\end{proof}

Without overparameterization, an RIP constant of $\delta<1/2$ is
both necessary and sufficient for global recovery. It is in fact possible
to prove global guarantees for RIP with all values of $\delta<1$
by incorporating the overparameterization ratio $r/r^{\star}$ as
a ``problem structure''. However, from the necessary condition,
we see that the rank parameter of the RIP assumption must also increase.
Here, we mention that rank overparameterization has only been recently
studied for matrix sensing~\citep{li2018algorithmic,zhuo2021computational};
its impact on the optimization landscape had been previously unknown.

\textbf{}

\subsubsection{Low-rank matrix recovery}

Our work is more distantly related to the more general problem of
\emph{low-rank matrix recovery}, of which matrix sensing is a particular
variant; see e.g. \citep{chi2019nonconvex} and the references therein.
The more general problem seeks to recover an unknown $M^{\star}\succeq0$
of low-rank $r^{\star}\ll n$ by minimizing an objective like $f(X)=\frac{1}{2}\|\AA(XX^{T})-b\|^{2}$
given possibly noisy measurements $b=\AA(M^{\star})+\varepsilon$,
but the choice of measurement operator $\AA$ usually does not satisfy
$(\delta,r+r^{\star})$-RIP. For certain instances like matrix completion
and robust PCA~\citep{ge2017nospurious}, the standard approach is
to design a regularized objective $f_{R}(X)=f(X)+R(X)$ whose global
minimizer $X_{R}^{\star}=\arg\min f_{R}(X)$ is guaranteed to recover
$M^{\star}$ up to a statistical error bound. In theory, if the regularized
objective can be written as $f_{R}(X)=\phi_{R}(XX^{T})$ with respect
to some $\phi_{R}$ that satisfies $(\delta,r+r^{\star})$-RIP, then
our results in this paper could be applied to guarantee global optimality,
and hence global recovery. Other instances like phase retrieval are
designed to work without a regularizer; without $(\delta,r+r^{\star})$-RIP,
our results are not applicable to these cases.

Recent works have studied overparameterization in the context of \emph{robust}
low-rank matrix recovery~\citep{ding2021rank,ma2022global}, wherein
the loss function is typically a nonsmooth $\ell_{1}$-based norm
like $f(X)=\|\AA(XX^{T})-b\|_{\ell_{1}}=\sum_{i}|\inner{A_{i}}{XX^{T}}-b_{i}|$.
A reviewer points out that, in the nonsmooth setting, a standard measure
of conditioning is the ratio of the Lipschitz constant over the modulus
of sharp growth. It is an important future direction to see if the
proof technique presented in this paper could be extended to the nonsmooth
setting.

\section{Notations}

\textbf{Basic linear algebra.} Lower-case letters are vectors and
upper-case letters are matrices. We use ``MATLAB notation'' in concatenating
vectors and matrices:

\[
[a,b]=\begin{bmatrix}a & b\end{bmatrix},\qquad[a;b]=\begin{bmatrix}a\\
b
\end{bmatrix},\qquad\diag(a,b)=\begin{bmatrix}a & 0\\
0 & b
\end{bmatrix},
\]
and the following short-hand to construct them:
\[
[x_{i}]_{i=1}^{n}=\begin{bmatrix}x_{1}\\
\vdots\\
x_{n}
\end{bmatrix},\qquad[x_{i,j}]_{i,j=1}^{m,n}=\begin{bmatrix}x_{1,1} & \cdots & x_{1,n}\\
\vdots & \ddots & \vdots\\
x_{m,1} & \cdots & x_{m,n}
\end{bmatrix}.
\]
Denote $\one=[1,1,\dots,1]^{T}$ as the vector-of-ones and $I=\diag(\one)$
as the identity matrix; we will infer their dimensions from context.
Denote $\langle X,Y\rangle=\tr(X^{T}Y)$ and $\|X\|_{F}^{2}=\langle X,X\rangle$
as the Frobenius inner product and norm. Denote $\nnz(x)$ as the
number of nonzero elements in $x$. The sets $\S^{n}\supset\S_{+}^{n}$
are the $n\times n$ real symmetric matrices, and the corresponding
positive semidefinite cone.

\textbf{Positive cones and projections.} The sets $\R^{n}\supset\R_{+}^{n}$
are the $n$ vectors, and the corresponding positive orthant. Denote
$x_{+}=\max\{0,x\}=\arg\min\{\|x-y\|:y\in\R_{+}^{n}\}$ as the projection
onto the positive orthant. 

\textbf{Vectorization and Kronecker product.} Denote $\vect(X)$ as
the usual column-stacking vectorizating operator, and $\mat(x)$ as
its inverse (the dimensions of the matrix are inferred from context).
The Kronecker product $\otimes$ is defined to satisfy the identity
$\vect(AXB^{T})=(B\otimes A)\vect(X)$.

\textbf{Pseudoinverse.} Denote the (Moore--Penrose) pseudoinverse
$A^{\dagger}\eqdef VS^{-1}U^{T}$ of matrix $A=USV^{T}$ where $U^{T}U=V^{T}V=I_{r}$
and $S\succ0$. Define $0^{\dagger}\eqdef0$.

\textbf{Error vector $\e$ and Jacobian matrix $\J_{X}$}. Given fixed
$X\in\R^{n\times r}$ and $Z\in\R^{n\times r}$, we denote $\e=\vect(XX^{T}-ZZ^{T})$
and implicitly define $\J_{X}\in\R^{n^{2}\times nr}$ to satisfy $\J_{X}\vect(V)=\vect(XV^{T}+VX^{T})$.
Note that the associate adjoint operator satisfies $\J_{X}^{T}\vect(M)=\vect((M+M^{T})X)$. 

\section{\label{sec:proof}Proof of the main result}

We will establish ``no spurious local minima'' guarantees by showing
that there cannot exist a counterexample that is simultaneously highly
overparameterized, with a large value of $r/r^{\star}$, and also
well-conditioned, with a small value of $\kappa$. 
\begin{definition}
The function $\phi:\S^{n}\to\R$ is said to be a $(\kappa,r,r^{\star})$-\emph{counterexample}
if it satisfies the following:
\begin{itemize}
\item The function $\phi$ is twice-differentiable, $L$-smooth, $\mu$-strongly
convex, and $L/\mu\le\kappa.$
\item The minimizer $M^{\star}=\arg\min_{M\succeq0}\phi(M)$ satisfies $\rank(M^{\star})=r^{\star}\le r$.
\item There exists an $n\times r$ spurious point $X$ with $XX^{T}\ne M^{\star}$
that satisfies $\nabla f(X)=0$ and $\nabla^{2}f(X)\succeq0$ for
$f(U)\eqdef\phi(UU^{T})$. 
\end{itemize}
\end{definition}
For a fixed $r$ and $r^{\star}$, our central claim in \thmref{main}
is that the condition number $\kappa^{\star}$ of the \emph{best-conditioned}
counterexample is exactly determined as 
\begin{equation}
\kappa^{\star}=\inf_{\phi:\S^{n}\to\R}\{\kappa:\phi\text{ is }(\kappa,r,r^{\star})\text{-counterexample}\}=1+2\sqrt{r/r^{\star}}.\label{eq:inf_kappa}
\end{equation}
The sufficient condition follows because there exists no $(\kappa,r,r^{\star})$-counterexample
with $\kappa<1+2\sqrt{r/r^{\star}}$, and so for a fixed $\kappa$,
we can always overparameterize $r/r^{\star}>\frac{1}{4}(\kappa-1)^{2}$
to eliminate all counterexamples. The neccessary condition holds because
a best-conditioned $(\kappa^{\star},r,r^{\star})$-counterexample
does exist.

To prove that $\kappa^{\star}$ is equal to the value stated in (\ref{eq:inf_kappa}),
we proceed by reformulating (\ref{eq:inf_kappa}) into two-stage optimization
problem 
\begin{align*}
\kappa^{\star} & =\inf_{X,Z\in\R^{n\times r}}\{\kappa(X,Z):\rank(Z)=r^{\star},\;XX^{T}\ne ZZ^{T}\}
\end{align*}
in which $\kappa(X,Z)$ corresponds to the best-conditioned counterexample
associated with a \emph{fixed} minimizer $ZZ^{T}$ and a \emph{fixed
}spurious point $X$:
\begin{equation}
\kappa(X,Z)=\inf_{\begin{subarray}{c}
\phi:\S^{n}\to\R\end{subarray}}\left\{ \kappa:\begin{array}{c}
\phi\text{ is }1\text{-strongly convex and }\kappa\text{-smooth,}\\
ZZ^{T}=\arg\min_{M\succeq0}\phi(M),\\
\nabla f(X)=0,\;\nabla^{2}f(X)\succeq0\text{ where }f(U)\eqdef\phi(UU^{T}).
\end{array}\right\} \label{eq:mu_reform}
\end{equation}
This two-stage reformulation is inspired by \citet{zhang2018much,zhang2019sharp},
who showed for a fixed $X,Z$ that finding the best-conditioned \emph{quadratic}
counterexample $\phi$ that satisfies $\nabla\phi(ZZ^{T})=0$ can
be exactly posed as a standard-form \emph{semidefinite program} or
SDP. Our proof extends this SDP formulation to a nonquadratic $\phi$
with a possibly nonzero $\nabla\phi(ZZ^{T})\succeq0$ at optimality
by introducing a slack variable.   

\begin{lemma}[SDP formulation]
\label{lem:sdp_formu}Define $\kappa(X,Z)$ as in (\ref{eq:mu_reform}).
We have $\kappa_{\ub}(X,Z)\ge\kappa(X,Z)\ge\kappa_{\lb}(X,Z)$ where
\begin{gather*}
\kappa_{\ub}(X,Z)\eqdef\min_{\kappa,\H}\left\{ \kappa:I\preceq\H\preceq\kappa I,\;\J_{X}^{T}\H\e=0,\;-2I_{r}\otimes\mat(\H\e)\preceq\J_{X}^{T}\H\J_{X}\right\} \\
\kappa_{\lb}(X,Z)\eqdef\min_{\kappa,s,\H}\left\{ \kappa:\begin{array}{c}
I\preceq\H\preceq\kappa I,\quad\mat(s)\succeq0,\quad\J_{Z}^{T}s=0,\\
\J_{X}^{T}(\H\e+s)=0,\quad-2I_{r}\otimes\mat(\H\e+s)\preceq\kappa\J_{X}^{T}\J_{X}
\end{array}\right\} 
\end{gather*}
and $\e=\vect(XX^{T}-ZZ^{T})$ and $\J_{U}\vect(V)=\vect(UV^{T}+VU^{T})$
for all $V$.
\end{lemma}
\begin{proof}
Denote $S_{M}$ as the convex gradient $\nabla\phi(M)$ evaluated
at $M$, and denote $\H_{M}$ as the matrix representation of the
convex Hessian operator $\nabla^{2}\phi(M)$ evaluated at $M$. It
follows from $\inner{\nabla^{2}\phi(M)[E]}E=\vect(E)^{T}\H_{M}\vect(E)$
that: 
\begin{equation}
\phi\text{ is }1\text{-strongly convex and }\kappa\text{-smooth}\iff I\preceq\H_{M}\preceq\kappa I\quad\text{for all }M\in\S^{n}.\label{eq:cond}
\end{equation}
Additionally, the Karush--Kuhn--Tucker conditions at $M^{\star}=ZZ^{T}$
read:
\begin{equation}
ZZ^{T}=\arg\min_{M\succeq0}\phi(M)\quad\iff\quad S_{ZZ^{T}}\succeq0,\quad S_{ZZ^{T}}Z=0.\label{eq:kktZZ}
\end{equation}
Note that strong duality holds in the above because $ZZ^{T}+\epsilon I$
is a strictly feasible point with a bounded objective due to $\kappa$-smoothness.
Next, we observe that the two gradient evaluations $S_{XX^{T}}$ and
$S_{ZZ^{T}}$ are related via the fundamental theorem of calculus:
\begin{gather*}
\nabla\phi(XX^{T})=\nabla\phi(ZZ^{T})+\int_{0}^{1}\nabla^{2}\phi((1-t)ZZ^{T}+tXX^{T})[XX^{T}-ZZ^{T}]\,\d t\\
\iff\quad\vect(S_{XX^{T}})=\vect(S_{ZZ^{T}})+\left[\int_{0}^{1}\H_{M(t)}\,\d t\right]\vect(XX^{T}-ZZ^{T}).
\end{gather*}
Denoting $s=\vect(S_{ZZ^{T}})$ and $\H_{\av}=\int_{0}^{1}\H_{M(t)}\,\d t$
and $\e=\vect(XX^{T}-ZZ^{T})$ allows us to rewrite the directional
derivatives of $f(U)\eqdef\phi(UU^{T})$ as follows: 
\begin{align}
\inner{\nabla f(X)}V & =\vect(V)^{T}\J_{X}^{T}(\H_{\av}\e+s),\nonumber \\
\inner{\nabla^{2}f(X)[V]}V & =\vect(V)^{T}[2I_{r}\otimes\mat(\H_{\av}\e+s)+\J_{X}^{T}\H_{XX^{T}}\J_{X}]\,\vect(V).\label{eq:vecdf}
\end{align}
Substituting (\ref{eq:cond}), (\ref{eq:vecdf}), and (\ref{eq:kktZZ})
into (\ref{eq:mu_reform}) yields exactly
\begin{equation}
\kappa(X,Z)=\inf_{\kappa,\phi:\S^{n}\to\R}\left\{ \kappa:\begin{array}{c}
I\preceq\H_{M}\preceq\kappa I\quad\text{for all }M\in\S^{n},\\
\mat(s)\succeq0,\quad\J_{Z}^{T}s=0,\quad\J_{X}^{T}(\H_{\av}\e+s)=0,\\
-2I_{r}\otimes\mat(\H_{\av}\e+s)\preceq\J_{X}^{T}\H_{XX^{T}}\J_{X}.
\end{array}\right\} \label{eq:kappa_sdp1}
\end{equation}
We obtain the upper-bound $\kappa_{\ub}(X,Z)$ from $\kappa(X,Z)$
in  (\ref{eq:kappa_sdp1}) by fixing $\H_{M}\equiv\H$ for all $M$,
and $S_{ZZ^{T}}\equiv\mat(s)=0$. We obtain the lower-bound $\kappa_{\lb}(X,Z)$
from $\kappa(X,Z)$ in (\ref{eq:kappa_sdp1}) by substituting the
relaxation 
\[
-2I_{r}\otimes\mat(\H_{\av}\e+s)\quad\preceq\quad\J_{X}^{T}\H_{XX^{T}}\J_{X}\quad\preceq\quad\kappa\cdot\J_{X}^{T}\J_{X},
\]
in order to avoid explicitly optimizing over $\H_{XX^{T}}$:
\begin{align}
\kappa(X,Z) & \ge\inf_{\kappa,\phi:\S^{n}\to\R}\left\{ \kappa:\begin{array}{c}
I\preceq\H_{M}\preceq\kappa I\quad\text{for all }M\in\S^{n},\\
\mat(s)\succeq0,\quad\J_{Z}^{T}s=0,\quad\J_{X}^{T}(\H_{\av}\e+s)=0,\\
-2I_{r}\otimes\mat(\H_{\av}\e+s)\preceq\kappa\cdot\J_{X}^{T}\J_{X}.
\end{array}\right\} \label{eq:kappasdp2}
\end{align}
Finally, we impose $I\preceq\H_{M}\preceq\kappa I$ over the average
Hessian $\H_{\av}=\int_{0}^{1}\H_{M(t)}\,\d t$ and then relax this
constraint over $\H_{M}$ for all other values of $M$.\qed
\end{proof}

\begin{remark}
We derived $\kappa_{\ub}(X,Z)$ by restricting the candidate counterexamples
$\phi$ in (\ref{eq:mu_reform}) to quadratics like $\phi(M)=\frac{1}{2}\vect(M-M^{\star})^{T}\H\vect(M-M^{\star})$.
Note that we also set $\nabla\phi(M^{\star})=0$ to avoid having to
optimize over it. \qed
\end{remark}
\begin{remark}
We derived $\kappa_{\lb}(X,Z)$ by substituting the following relaxation,
which is implied by the $L$-smoothness of $\phi$, into (\ref{eq:mu_reform}):
\[
\inner{\nabla^{2}f(X)[V]}V\ge0\implies2\inner{\nabla\phi(XX^{T})}{VV^{T}}+L\cdot\|XV^{T}+VX^{T}\|_{F}^{2}\ge0.
\]
After making this relaxation, the best-conditioned $\phi$ turns out
to be a quadratic like $\phi(M)=s^{T}\vect(M-M^{\star})+\frac{1}{2}\vect(M-M^{\star})^{T}\H\vect(M-M^{\star}).$
However, $\phi$ might not be a valid counterexample; due to the relaxation,
it is not necessarily feasible for (\ref{eq:mu_reform}). \qed
\end{remark}
We prove $\kappa^{\star}\le1+2\sqrt{r/r^{\star}}$ by stating an explicit
choice of $X$ and $Z$ that heuristically minimizes the upper-bound
$\kappa_{\ub}(X,Z)$ in \lemref{sdp_formu} over $n\times r$ matrices
$X$ and full-rank $n\times r^{\star}$ matrices $Z$. We defer its
proof to \secref{counterexamples}. 
\begin{lemma}[Heuristic upper-bound]
\label{lem:heuristic}Let $[Q_{1},Q_{2}]$ have orthonormal columns
with $Q_{1}\in\R^{n\times r}$ and $Q_{2}\in\R^{n\times r^{\star}}$.
Then, we have
\[
\kappa_{\ub}(X,Z)\le1+2\sqrt{r/r^{\star}}\quad\text{where }X=Q_{1},\quad Z=\sqrt{1+\sqrt{r/r^{\star}}}Q_{2}.
\]
\end{lemma}
We prove $\kappa^{\star}\ge1+2\sqrt{r/r^{\star}}$ by deriving a closed-form
lower-bound on $\kappa_{\lb}(X,Z)$ in \lemref{sdp_formu}, and then
analytically minimizing it over $n\times r$ matrices $X$ and full-rank
$n\times r^{\star}$ matrices $Z$.  We defer the derivation of this
lower-bound to \secref{simple}. 

\begin{lemma}[Closed-form lower-bound]
\label{lem:abdef}For $X,Z\in\R^{n\times r}$ such that $XX^{T}\ne ZZ^{T}$,
let $Z_{\perp}=(I-XX^{\dagger})Z$. We have
\begin{gather*}
\kappa_{\lb}(X,Z)\ge\begin{cases}
{\displaystyle \frac{1+\sqrt{1-\alpha^{2}}}{1-\sqrt{1-\alpha^{2}}}} & \text{if }\beta\ge{\displaystyle \frac{\alpha}{1+\sqrt{1-\alpha^{2}}}},\\
{\displaystyle \frac{1-\alpha\beta}{(\alpha-\beta)\beta}} & \text{if }\beta\le{\displaystyle \frac{\alpha}{1+\sqrt{1-\alpha^{2}}}},
\end{cases}\\
\text{ where }\qquad\alpha=\frac{\|Z_{\perp}Z_{\perp}^{T}\|_{F}}{\|XX^{T}-ZZ^{T}\|_{F}},\qquad\beta=\frac{\lambda_{\min}(X^{T}X)}{\|XX^{T}-ZZ^{T}\|_{F}}{\displaystyle \frac{\tr(Z_{\perp}Z_{\perp}^{T})}{\|Z_{\perp}Z_{\perp}^{T}\|_{F}}}.
\end{gather*}
\end{lemma}
\begin{remark}
The definition of $\beta$ becomes ambiguous when $Z_{\perp}=0$.
This is without loss of precision, because $Z_{\perp}=0$ implies
$\alpha=0$, and therefore $\kappa_{\lb}(X,Z)=1$ regardless of the
value of $\beta$. 
\end{remark}
Our main difficulty in our proof is the need to minimize the closed-form
lower-bound in \lemref{abdef} over all possible choices of $X$ and
$Z$. In the rank-1 case, this easily follows by substituting the
symmetry invariants $\rho=\|x\|/\|z\|$ and $\varphi=\arccos(\frac{x^{T}z}{\|x\|\|z\|})$,
as in 
\[
\alpha=\frac{\sin^{2}\varphi}{\sqrt{(1-\rho^{2})^{2}+2\rho^{2}\sin^{2}\varphi}},\qquad\beta=\frac{\rho^{2}}{\sqrt{(1-\rho^{2})^{2}+2\rho^{2}\sin^{2}\varphi}},
\]
and then explicitly minimizing the lower-bound with respect to $\rho$
and $\varphi$; see \citet[Theorem~3]{zhang2019sharp}. In the possibly
overparameterized rank-$r$ case, however, the same approach would
force us to solve a nonconvex optimization problem over up to $rr^{\star}+r+r^{\star}-1$
symmetry invariants. Moreover, it is unclear how we can preserve the
nonconvex rank equality constraint $\rank(Z)=r^{\star}$ once the
problem is reformulated into symmetry invariants. 

The main innovation in our proof is to relax the explicit dependence
of $\alpha$ and $\beta$ on their arguments $X,Z$, and to minimize
the lower-bound in \lemref{abdef} directly over $\alpha,\beta$ as
variables. Adopting a classic strategy from integer programming, we
tighten this relaxation by introducing a \emph{valid inequality}:
an inequality constraint on $\alpha,\beta$ that remains valid for
all choices of $X,Z\in\R^{n\times r}$ with $\rank(Z)=r^{\star}$.
Below, we write $(\cdot)_{+}=\max\{\cdot,0\}$.
\begin{lemma}[Valid inequality]
\label{lem:ab1}For all $X,Z\in\R^{n\times r}$ with $\rank(Z)=r^{\star}$,
the parameters $\alpha$ and $\beta$ defined in \Lemref{abdef} satisfy
$\alpha^{2}+(r/r^{\star})\beta^{2}\le1+[(\beta-\alpha)_{+}]^{2}.$
\end{lemma}
We will soon prove \lemref{ab1} in \secref{FeasRegion} below. Surprisingly,
this very simple valid inequality is all that is needed to minimize
our simple lower-bound $\kappa_{\lb}(X,Z)$ to global optimality.
The proof of \thmref{main} quickly follows from \lemref{heuristic},
\lemref{abdef}, and \lemref{ab1}. 
\begin{proof}
[Theorem~\ref{thm:main}]For a fixed $r$ and $r^{\star}$, recall
that we have defined $\kappa^{\star}$ as the condition number of
the \emph{best-conditioned} counterexample, and formulated it as the
optimal value for a two-stage optimization problem of the following
form
\begin{align*}
\kappa^{\star} & =\inf_{\phi:\S^{n}\to\R}\{\kappa:\phi\text{ is }(\kappa,r,r^{\star})\text{-counterexample}\}\\
 & =\inf_{X,Z\in\R^{n\times r}}\{\kappa(X,Z):\rank(Z)=r^{\star},\;XX^{T}\ne ZZ^{T}\}
\end{align*}
We claim that the minimum of the following expression over $\alpha,\beta\ge0$
\[
\gamma(\alpha,\beta)\eqdef\begin{cases}
\frac{1+\sqrt{1-\alpha^{2}}}{1-\sqrt{1-\alpha^{2}}} & \text{if }\beta\ge{\displaystyle \frac{\alpha}{1+\sqrt{1-\alpha^{2}}}},\\
{\displaystyle \frac{1-\alpha\beta}{(\alpha-\beta)\beta}} & \text{if }\beta\le{\displaystyle \frac{\alpha}{1+\sqrt{1-\alpha^{2}}}},
\end{cases}
\]
subject to the valid inequality from \lemref{ab1} is given:
\begin{equation}
1+2\sqrt{r/r^{\star}}=\min_{\alpha,\beta\ge0}\{\gamma(\alpha,\beta):\alpha^{2}+(r/r^{\star})\min\{\alpha^{2},\beta^{2}\}\le1\}.\label{eq:minab}
\end{equation}
Taking (\ref{eq:minab}) to be true, substituting into \lemref{abdef}
and evoking \lemref{heuristic} yields
\[
1+2\sqrt{r/r^{\star}}\le\inf_{\begin{subarray}{c}
XX^{T}\ne ZZ^{T}\\
\rank(Z)=r^{\star}
\end{subarray}}\kappa_{\lb}(X,Z)\le\kappa^{\star}\le\kappa_{\ub}(X,Z)\le1+2\sqrt{r/r^{\star}}.
\]
The upper-bound proves necessity: if $r/r^{\star}\le\frac{1}{4}(\kappa-1)^{2}$,
then there exists a $(\kappa,r,r^{\star})$-counterexample. The lower-bound
proves sufficiency: if $r>\frac{1}{4}(\kappa-1)^{2}\,r^{\star}$ and
$r^{\star}\le r<n$, then there exists no $(\kappa,r,r^{\star})$-counterexample
to refute a ``no spurious local minima'' guarantee.

We now prove the claim in (\ref{eq:minab}). We first minimize $\gamma(\alpha,\beta)$
subject to $\beta\ge\frac{\alpha}{1+\sqrt{1-\alpha^{2}}}$ and find
via monotonicity with respect to $\alpha$ that
\[
\min_{0\le\alpha\le1}\left\{ \frac{1+\sqrt{1-\alpha^{2}}}{1-\sqrt{1-\alpha^{2}}}:\beta\ge{\displaystyle \frac{\alpha}{1+\sqrt{1-\alpha^{2}}}}\right\} ={\displaystyle \frac{1-\alpha\beta}{(\alpha-\beta)\beta}}.
\]
Indeed, the minimizer $\alpha$ lies at the boundary $\beta=\frac{\alpha}{1+\sqrt{1-\alpha^{2}}}=\frac{1-\sqrt{1-\alpha^{2}}}{\alpha}$,
where
\[
\frac{1+\sqrt{1-\alpha^{2}}}{1-\sqrt{1-\alpha^{2}}}=\frac{\sqrt{1-\alpha^{2}}+(1-\alpha^{2})}{\sqrt{1-\alpha^{2}}-(1-\alpha^{2})}=\frac{1+\sqrt{1-\alpha^{2}}-\alpha^{2}}{\alpha^{2}-(1-\sqrt{1-\alpha^{2}})}={\displaystyle \frac{\beta^{-1}-\alpha}{\alpha-\beta}.}
\]
Next, we minimize $\gamma(\alpha,\beta)$ subject to $\beta\le\frac{\alpha}{1+\sqrt{1-\alpha^{2}}}\le\alpha$
and the valid inequality:
\begin{align*}
 & \min_{0\le\beta\le\alpha}\left\{ {\displaystyle \frac{1-\alpha\beta}{(\alpha-\beta)\beta}}:\alpha^{2}+\frac{r}{r^{\star}}\beta^{2}\le1\right\} \\
= & \min_{\beta\ge0,\;\rho\ge1}\left\{ {\displaystyle \frac{1-\rho\beta^{2}}{(\rho-1)\beta^{2}}}:(\rho^{2}+\frac{r}{r^{\star}})\beta^{2}\le1\right\} =\min_{\rho\ge1}\left\{ \rho+\frac{r}{r^{\star}}\frac{1}{\rho-1}\right\} =1+2\sqrt{\frac{r}{r^{\star}}}.
\end{align*}
We obtain the second line from the first by substituting $\alpha=\rho\beta$
for $\rho\ge1$. The minimizers are $\beta^{2}=(\rho^{2}+\frac{r}{r^{\star}})^{-1}$
and $\rho=1+\sqrt{r/r^{\star}}$. \qed
\end{proof}

\section{\label{sec:FeasRegion}Proof of the valid inequality over $\alpha$
and $\beta$ (\lemref{ab1})}

The main innovation in our proof is \lemref{ab1}, which provided
a valid inequality that allowed us to relax the dependence of the
parameters $\alpha$ and $\beta$ in our closed-form lower-bound (\lemref{abdef})
with respect to $X,Z$, while simultaneously capturing the rank \emph{equality}
constraint $r^{\star}=\rank(Z)$. The lemma in turn crucially depends
on the following generalization of the classic result of \citet{eckart1936approximation}.
\begin{theorem}[Regularized Eckart--Young]
\label{thm:EckartYoung}Given $A\in\S_{+}^{n}$ and $B\in\S_{+}^{r}$
with $r\le n$, let $A=\sum_{i=1}^{n}s_{i}u_{i}u_{i}^{T}$ and $B=\sum_{i=1}^{r}d_{i}v_{i}v_{i}^{T}$
denote the usual orthonormal eigendecompositions with $s_{1}\ge\cdots\ge s_{n}\ge0$
and $0\le d_{1}\le\cdots\le d_{r}.$ Then, 
\begin{align}
\min_{Y\in\R^{n\times r}}\left\{ \|A-YY^{T}\|_{F}^{2}+2\langle B,Y^{T}Y\rangle\right\}  & =\sum_{i=1}^{n}s_{i}^{2}-\sum_{i=1}^{r}[(s_{i}-d_{i})_{+}]^{2}\label{eq:ey-sol}
\end{align}
with minimizer $Y^{\star}=\sum_{i=1}^{r}u_{i}v_{i}^{T}\sqrt{(s_{i}-d_{i})_{+}}$
where $(\cdot)_{+}\eqdef\max\{0,\cdot\}$. 
\end{theorem}
Setting $B=0$ in Theorem~\ref{thm:EckartYoung} recovers the original
Eckart--Young Theorem: The best rank-$r$ approximation $YY^{T}\approx A$
in Frobenius norm is the truncated singular value decomposition $YY^{T}=\sum_{i=1}^{r}s_{i}u_{i}u_{i}^{T},$
with approximation error $\|A-YY^{T}\|_{F}^{2}=\sum_{i=r+1}^{n}s_{i}^{2}$.
Hence, $B\ne0$ may be viewed as a regularizer that prevents $Y^{\star}$
from becoming excessively large. 

Our proof of \thmref{EckartYoung} follows a trace inequality of Ruhe~\citep{ruhe1970perturbation}
(see also \citep[Chapter~9]{marshall2011inequalities}).\footnote{This simplified proof via Ruhe's inequality is due to an anonymous
reviewer.} The inequality states that, for $n\times n$ positive semidefinite
matrices $A,B$, we have
\begin{equation}
\sum_{i=1}^{n}\lambda_{i}(A)\lambda_{n-i+1}(B)\le\inner AB\le\sum_{i=1}^{n}\lambda_{i}(A)\lambda_{i}(B)\label{eq:Ruhe}
\end{equation}
where the eigenvalues are ordered $\lambda_{1}(A)\ge\cdots\ge\lambda_{n}(A)$
and $\lambda_{1}(B)\ge\cdots\ge\lambda_{n}(B)$. 
\begin{proof}
[\thmref{EckartYoung}]Recall that $A=\sum_{i=1}^{n}s_{i}u_{i}u_{i}^{T}$
and $B=\sum_{i=1}^{r}d_{i}v_{i}v_{i}^{T}$ with $s_{1}\ge\cdots\ge s_{n}\ge0$
and $0\le d_{1}\le\cdots\le d_{r}.$ It follows from (\ref{eq:Ruhe})
that
\begin{gather*}
\inner A{YY^{T}}\le\sum_{i=1}^{n}\lambda_{i}(YY^{T})\cdot\lambda_{i}(A)=\sum_{i=1}^{r}s_{i}\sigma_{i}^{2}(Y),\\
\inner B{Y^{T}Y}\ge\sum_{i=1}^{r}\lambda_{i}(Y^{T}Y)\cdot\lambda_{n-i+1}(B)=\sum_{i=1}^{r}d_{i}\sigma_{i}^{2}(Y).
\end{gather*}
Substituting yields the following lower-bound, which is valid for
all choices of $Y$:
\begin{align*}
\|A-YY^{T}\|_{F}^{2}+2\inner B{Y^{T}Y}= & \|A\|_{F}^{2}+\|YY^{T}\|_{F}^{2}-2\inner A{YY^{T}}+2\inner B{Y^{T}Y}\\
\ge & \sum_{i=1}^{n}s_{i}^{2}+\sum_{i=1}^{r}\sigma_{i}^{4}(Y)-2\sum_{i=1}^{r}s_{i}\sigma_{i}^{2}(Y)+2\sum_{i=1}^{r}d_{i}\sigma_{i}^{2}(Y)\\
\ge & \sum_{i=1}^{n}s_{i}^{2}-\sum_{i=1}^{r}[(s_{i}-d_{i})_{+}]^{2}.
\end{align*}
The final line applied the following result as a lower-bound:
\begin{align*}
\sigma_{i}^{4}(Y)-2(s_{i}-d_{i})\sigma_{i}^{2}(Y)\ge & \min_{x_{i}\ge0}x_{i}^{2}-2(s_{i}-d_{i})x_{i}\\
= & \min_{x_{i}\ge0}[x_{i}-(s_{i}-d_{i})]^{2}-(s_{i}-d_{i})^{2}\\
= & [(s_{i}-d_{i})_{+}-(s_{i}-d_{i})]^{2}-(s_{i}-d_{i})^{2}=-[(s_{i}-d_{i})_{+}]^{2}.
\end{align*}
Finally, we verify that $Y^{\star}=\sum_{i=1}^{n}u_{i}v_{i}\sqrt{(s_{i}-d_{i})_{+}}$
attains the lower-bound derived above. \qed
\end{proof}

The proof of Lemma~\ref{lem:ab1} quickly follows from Theorem~\ref{thm:EckartYoung}.
We will also need the following two technical lemmas.
\begin{lemma}
\label{lem:numcard}If $x\ge0$ and $\one^{T}x\le\|x\|^{2}$, then
$\one^{T}\left(I-xx^{T}/\|x\|^{2}\right)\one\ge\|(\one-x)_{+}\|^{2}$.
\end{lemma}
\begin{proof}
Define $u=\alpha x,$ where $\alpha=\one^{T}x/\|x\|^{2}$ is chosen
so that $\mathbf{1}^{T}u=\|u\|^{2}$. First, verify that $\mathbf{1}^{T}(I-xx^{T}/\|x\|^{2})\mathbf{1}=\|\mathbf{1}-u\|^{2}\ge\|(\mathbf{1}-u)_{+}\|^{2}$.
Next, observe that $\psi(\alpha)\equiv[(1-\alpha t)_{+}]^{2}$ is
a decreasing function of $\alpha\ge0$ when $t\ge0$, and therefore
$\|(\one-\alpha x)_{+}\|^{2}\ge\|(\one-x)_{+}\|^{2}$ for $x\ge0$
and $0\le\alpha\le1$. Finally, the hypotheses $x\ge0$ and $\one^{T}x\le\|x\|^{2}$
ensure that $0\le\alpha\le1$. \qed
\end{proof}

\begin{lemma}
\label{lem:keyclaim}Given $s,d\in\R_{+}^{n}$, let $s_{i}\ge s_{\lb}\ge0$
for all $i$. Then we have
\begin{align*}
\|s\|^{2}-\|(s-d)_{+}\|^{2}\ge & \begin{cases}
s_{\lb}^{2}\frac{(\one^{T}d)^{2}}{\|d\|^{2}} & \text{if }s_{\lb}\one^{T}d\le\|d\|^{2},\\
2s_{\lb}\one^{T}d-\|d\|^{2} & \text{if }s_{\lb}\one^{T}d\ge\|d\|^{2}.
\end{cases}
\end{align*}
\end{lemma}
\begin{proof}
First, we observe that $\psi(t)\equiv t^{2}-[(t-\alpha)_{+}]^{2}$
is an increasing function of $t\ge0$. It follows from element-wise
monotonicity that
\[
\|s\|^{2}-\|(s-d)_{+}\|^{2}\ge\|s_{\lb}\one\|^{2}-\|(s_{\lb}\one-d)_{+}\|^{2}=s_{\lb}^{2}\left(\|\one\|^{2}-\|(\one-d/s_{\lb})_{+}\|^{2}\right)
\]
If $s_{\lb}\one^{T}d\le\|d\|^{2}$, then write $x=d/s_{\lb}\ge0$,
and observe that $\one^{T}x\le\|x\|^{2}$ holds. Therefore, applying
\lemref{numcard} yields
\begin{align*}
s_{\lb}^{2}\left(\|\one\|^{2}-\|(\one-d/s_{\lb})_{+}\|^{2}\right) & =s_{\lb}^{2}\left(\|\one\|^{2}-\|(\one-x)_{+}\|^{2}\right)\ge s_{\lb}^{2}\one^{T}\left(xx^{T}/\|x\|^{2}\right)\one\\
 & =s_{\lb}^{2}\left([\one^{T}(d/s_{\lb})]^{2}/\|d/s_{\lb}\|^{2}\right)=s_{\lb}^{2}(\one^{T}d)^{2}/\|d\|^{2}.
\end{align*}
If otherwise $s_{\lb}\one^{T}d\ge\|d\|^{2}$, then we have 
\begin{align*}
s_{\lb}^{2}\left(\|\one\|^{2}-\|(\one-d/s_{\lb})_{+}\|^{2}\right)\ge s_{\lb}^{2}\left(\|\one\|^{2}-\|\one-d/s_{\lb}\|^{2}\right)=2s_{\lb}\one^{T}d-\|d\|^{2}.
\end{align*}
\qed
\end{proof}

Instead of proving \lemref{ab1}, we will prove the following equivalent
statement.
\begin{lemma}
\label{lem:ab2}For $X,Z\in\R^{n\times r}$, define the parameters
\[
\alpha=\frac{\|Z_{\perp}Z_{\perp}^{T}\|_{F}}{\|XX^{T}-ZZ^{T}\|_{F}},\qquad\beta=\frac{\lambda_{\min}(X^{T}X)}{\|XX^{T}-ZZ^{T}\|_{F}}\cdot\frac{\tr(Z_{\perp}Z_{\perp}^{T})}{\|Z_{\perp}Z_{\perp}^{T}\|_{F}}.
\]
Then, $\alpha^{2}+(r/p)\beta^{2}\le1+[(\beta-\alpha)_{+}]^{2}$ holds
for all $p$ satisfying $\rank(Z)\le p\le r$.
\end{lemma}
Assuming the above to be true, setting $p=\rank(Z)$ yields exactly
\lemref{ab1}. Conversely, if the above is true for $p=\rank(Z)$,
then it is obviously true for any $p\ge\rank(Z)$, given that the
left-hand side decreases monotonically with increasing $p$.
\begin{proof}
The two parameters $\alpha$ and $\beta$ are unitarily invariant
with respect to $X,Z$, so we may assume without loss of generality
that $X,Z$ are partitioned as
\[
X=\begin{bmatrix}X_{1}\\
0
\end{bmatrix}\qquad Z=\begin{bmatrix}Z_{1}\\
Z_{2}
\end{bmatrix}\qquad\text{where }X_{1}\in\R^{r\times r},\;Z_{1}\in\R^{r\times p},\;Z_{2}\in\R^{(n-r)\times p}.
\]
(Otherwise, take the QR decomposition $X=QR$ and note that $\alpha,\beta$
remain unchanged with $X\gets Q^{T}X$ and $Z\gets Q^{T}Z$.) Denote
$s\in\R^{r}$ as the $r$ eigenvalues of $X_{1}X_{1}^{T}$ and $d\in\R^{p}$
as the $p$ eigenvalues of $Z_{2}^{T}Z_{2}$ 
\begin{align*}
s_{i}\eqdef & \lambda_{i}(X_{1}X_{1}^{T})=\lambda_{i}(XX^{T}),\quad s_{1}\ge s_{2}\ge\cdots\ge s_{r}\ge0,\\
d_{i}\eqdef & \lambda_{i}(Z_{2}^{T}Z_{2})=\lambda_{i}(Z_{\perp}Z_{\perp}^{T}),\quad0\le d_{1}\le d_{2}\le\cdots\le d_{p}.
\end{align*}
Writing $\e=\vect(XX^{T}-ZZ^{T})$ yields
\begin{align*}
\alpha & =\frac{\|Z_{2}Z_{2}^{T}\|_{F}}{\|XX^{T}-ZZ^{T}\|_{F}}=\frac{\|d\|_{F}}{\|\e\|_{F}}, & \beta & =\frac{\lambda_{\min}(X_{1}^{T}X_{1})}{\|XX^{T}-ZZ^{T}\|_{F}}\cdot\frac{\tr(Z_{2}Z_{2}^{T})}{\|Z_{2}Z_{2}^{T}\|_{F}}=\frac{s_{r}\cdot\one^{T}d}{\|\e\|\|d\|}.
\end{align*}
We expand $\|XX^{T}-ZZ^{T}\|_{F}^{2}$ block-wise and evoke our Regularized
Eckart--Young Theorem (\thmref{EckartYoung}) over the $r\times p$
matrix block $Z_{1}$:
\begin{align*}
\|XX^{T}-ZZ^{T}\|_{F}^{2} & =\|X_{1}X_{1}^{T}-Z_{1}Z_{1}^{T}\|_{F}^{2}+2\langle Z_{2}^{T}Z_{2},Z_{1}^{T}Z_{1}\rangle+\|Z_{2}Z_{2}^{T}\|^{2}\\
 & \ge\min_{Z_{1}\in\R^{r\times p}}\{\|X_{1}X_{1}^{T}-Z_{1}Z_{1}^{T}\|_{F}^{2}+2\langle Z_{2}^{T}Z_{2},Z_{1}^{T}Z_{1}\rangle\}+\|Z_{2}Z_{2}^{T}\|^{2}\\
 & =\|s\|^{2}-\|(s'-d)_{+}\|^{2}+\|d\|^{2}
\end{align*}
where $s'=(s_{1},s_{2},\dots,s_{p})$. Now, evoking our technical
lemma (\lemref{keyclaim}) with $s_{\lb}\equiv s_{r}$ yields the
following two lower-bounds on $\|XX^{T}-ZZ^{T}\|_{F}^{2}=\|\e\|^{2}$:\begin{subequations}\label{eq:key-claim}
\begin{gather}
s_{r}\one^{T}d\le\|d\|^{2}\quad\implies\quad\|\e\|^{2}\ge s_{r}^{2}\frac{(\one^{T}d)^{2}}{\|d\|^{2}}+(r-p)s_{r}^{2}+\|d\|^{2},\label{eq:key-claim1}\\
s_{r}\one^{T}d\ge\|d\|^{2}\quad\implies\quad\|\e\|^{2}\ge2s_{r}\one^{T}d-\|d\|^{2}+(r-p)s_{r}^{2}+\|d\|^{2}.\label{eq:key-claim2}
\end{gather}
\end{subequations}Substituting $s_{r}=\beta\cdot\|\e\|\|d\|/(\one^{T}d)$
yields
\[
(r-p)s_{r}^{2}=(r-p)\frac{\|d\|^{2}}{(\one^{T}d)^{2}}\beta^{2}\|\e\|^{2}\ge(r-p)\frac{1}{\nnz(d)}\beta^{2}\|\e\|^{2}\ge\frac{(r-p)}{p}\beta^{2}\|\e\|^{2},
\]
and substituting into (\ref{eq:key-claim}) yields\begin{subequations}
\begin{gather}
\beta\le\alpha\quad\implies\quad1\ge\beta^{2}+\alpha^{2}+\frac{(r-p)}{p}\beta^{2}=\alpha^{2}+\frac{r}{p}\beta^{2},\label{eq:key-claim3}\\
\beta\ge\alpha\quad\implies\quad1\ge2\alpha\beta+\frac{(r-p)}{p}\beta^{2}=\alpha^{2}+\frac{r}{p}\beta^{2}-(\beta-\alpha)^{2},\label{eq:key-claim4}
\end{gather}
\end{subequations}as desired. \qed
\end{proof}

\section{\label{sec:simple}Proof of the closed-form lower-bound (\lemref{abdef})}

Given fixed $X,Z\in\R^{n\times r}$, define the error vector $\e=\vect(XX^{T}-ZZ^{T})$.
Recall that we define $\J_{X}$ with respect to a given $X\in\R^{n\times r}$
to implicitly satisfy 
\begin{equation}
\J_{X}\vect(V)=\vect(XV^{T}+VX^{T}),\qquad\J_{X}^{T}\vect(M)=\vect((M+M^{T})X)\label{eq:Jxdef}
\end{equation}
for all $V\in\R^{n\times r}$ and $M\in\R^{n\times n}$. 
\begin{lemma}
\label{lem:I-JJ+}For any $X\in\R^{n\times r}$, the matrix $\J_{X}$
defined in (\ref{eq:Jxdef}) satisfies $I-\J_{X}\J_{X}^{\dagger}=(I-XX^{\dagger})\otimes(I-XX^{\dagger})$
where $\dagger$ denotes the pseudoinverse. 
\end{lemma}
\begin{proof}
Let $X=QR$ be the usual QR decomposition with orthogonal complement
$Q_{\perp}$, so that $XX^{\dagger}=QQ^{T}$ and $I-XX^{\dagger}=Q_{\perp}Q_{\perp}^{T}$.
Observe that for any $M\in\R^{n\times n}$:
\[
\J_{X}^{T}\vect(M)=0\iff(M+M^{T})Q=0\iff M=Q_{\perp}\hat{M}Q_{\perp}^{T}.
\]
Vectorizing the final equality as $\vect(M)=(Q_{\perp}\otimes Q_{\perp})\vect(\hat{M})$
yields 
\begin{align*}
(I-\J_{X}\J_{X}^{\dagger})b & =\arg\min_{v}\{\|v-b\|:\J_{X}^{T}v=0\}\\
 & =\arg\min_{v}\{\|v-b\|:v=(Q_{\perp}\otimes Q_{\perp})\hat{v}\}\\
 & =(Q_{\perp}\otimes Q_{\perp})(Q_{\perp}\otimes Q_{\perp})^{T}b=(Q_{\perp}Q_{\perp}^{T}\otimes Q_{\perp}Q_{\perp}^{T})b.\quad\qed
\end{align*}
\end{proof}

To lower-bound $\kappa_{\lb}(X,Z)$ as defined in \lemref{sdp_formu},
we will plug a carefully-chosen heuristic solution into its Lagrangian
dual. This dual problem is stated in the lemma below; see \appref{tradeoff}
for its derivation.
\begin{lemma}
\label{lem:tradeoff}We have
\[
\kappa_{\lb}(X,Z)=\max_{t\ge0}\frac{1+\cos\theta(t)}{2t+1-\cos\theta(t)}
\]
where $\cos\theta(t)$ is defined as follows
\begin{equation}
\cos\theta(t)=\max_{y,W_{i,j}}\left\{ \e^{T}\f:\begin{array}{c}
\f\equiv\J_{X}y-\sum_{i=1}^{r}\vect(W_{i,i}),\quad\|\e\|\|\f\|=1,\\
(I-ZZ^{\dagger})\,\mat(\f)\,(I-ZZ^{\dagger})\succeq0,\\
\langle\J_{X}^{T}\J_{X},W\rangle=2t,\quad W=[W_{i,j}]_{i,j=1}^{r}\succeq0,
\end{array}\right\} \label{eq:cosprob}
\end{equation}
with respect to $y\in\R^{nr}$ and $W_{i,j}\in\R^{n\times n}$ for
$i,j\in\{1,2,\dots,r\}$.
\end{lemma}
Our heuristic choice of $y$ and $W_{i,j}$ is motivated by the optimal
solution to the much simpler rank-1 case previously derived by \citet{zhang2019sharp}.
This choice is not necessarily optimal for the rank-$r$ case, but
it does remain simple enough for the final maximization over $t\ge0$
in \lemref{tradeoff} to be solved in closed-form. 

We begin by recalling the core insights used to derive the rank-1
solution. Our goal of maximizing $\e^{T}(\J_{X}y-w)$ where $w=\sum_{i=1}^{r}\vect(W_{i,i})$
can equivalently be understood as minimizing $\|\J_{X}y-w-\e\|$ over
$y$ and $w$, subject to $\mat(w)\succeq0$. This is a projection
problem, so a basic idea is to decompose $\e$ into orthogonal components:
\begin{equation}
\J_{X}y=\gamma_{1}\cdot(\J_{X}\J_{X}^{\dagger})\e,\qquad w=-\gamma_{2}\cdot(I-\J_{X}\J_{X}^{\dagger})\e.\label{eq:ysol}
\end{equation}
We verify that $\mat(w)\succeq0$ indeed holds for $\gamma_{2}\ge0$
via \lemref{I-JJ+}
\begin{equation}
\mat((I-\J_{X}\J_{X}^{\dagger})\e)=(I-XX^{\dagger})(XX^{T}-ZZ^{T})(I-XX^{\dagger})=-Z_{\perp}Z_{\perp}^{T}\label{eq:e2}
\end{equation}
where we have written $Z_{\perp}=(I-XX^{\dagger})Z$. 

Equation (\ref{eq:ysol}) specifies a unique choice of $y$, but leaves
several possible valid choices of $W_{i,j}$. We need the following
claim to make a concrete choice.
\begin{claim}
Let $W=\sum_{j=1}^{m}\vect(V_{j})\vect(V_{j})^{T}$ where $V_{1},\dots,V_{m}\in\R^{n\times r}$.
Partition $W$ into $r\times r$ blocks of $n\times n$, as in $W=[W_{i,j}]_{i,j=1}^{r}$
for $W_{i,j}\in\R^{n\times n}$. Then, 
\[
\sum_{i=1}^{r}W_{i,i}=\sum_{j=1}^{m}V_{j}V_{j}^{T},\qquad\inner{\J_{X}^{T}\J_{X}}W=\sum_{j=1}^{m}\|XV_{j}^{T}+V_{j}X^{T}\|_{F}^{2}.
\]
\end{claim}
Following the claim, the obvious first attempt is to pick $W=\gamma_{2}\cdot\vect(Z_{\perp})\vect(Z_{\perp})^{T}$.
Indeed, this yields $\sum_{i=1}^{r}W_{i,i}=\gamma_{2}\cdot Z_{\perp}Z_{\perp}^{T}$
as specified in (\ref{eq:ysol}), and is in fact optimal for the rank-1
case. In the rank-$r$ case, however, the choice of $W$ can be improved
by making $\inner{\J_{X}^{T}\J_{X}}W$ as small as possible. A more
subtle second attempt is to pick 
\begin{equation}
W=\gamma_{2}\cdot\sum_{j=1}^{r}\vect(Z_{\perp}e_{j}v_{r}^{T})\vect(Z_{\perp}e_{j}v_{r}^{T})^{T}\label{eq:Wsol}
\end{equation}
where $v_{r}$ is the $r$-th unit eigenvector of $X^{T}X$. Indeed,
we can verify that (\ref{eq:ysol}) is also satisfied:
\begin{gather}
\sum_{i=1}^{r}W_{i,i}=\gamma_{2}\cdot\sum_{j=1}^{r}Z_{\perp}e_{i}v_{r}^{T}v_{r}e_{i}^{T}Z_{\perp}^{T}=\gamma_{2}\cdot Z_{\perp}\left(\sum_{i=1}^{r}e_{i}e_{i}^{T}\right)Z_{\perp}^{T}=\gamma_{2}\cdot Z_{\perp}Z_{\perp}^{T}.\label{eq:Wsol1}
\end{gather}
Moreover, we have
\begin{align}
\langle\J_{X}^{T}\J_{X},W\rangle & =\gamma_{2}\cdot\sum_{i=1}^{m}\|\J_{X}\vect(Z_{\perp}e_{i}v_{r}^{T})\|^{2}=\gamma_{2}\cdot\sum_{i=1}^{m}\|Xv_{r}e_{i}^{T}Z_{\perp}^{T}+Z_{\perp}e_{i}v_{r}^{T}X^{T}\|_{F}^{2}\nonumber \\
 & =\gamma_{2}\cdot\sum_{i=1}^{m}2\|Z_{\perp}e_{i}v_{r}^{T}X^{T}\|_{F}^{2}=2\gamma_{2}\cdot\lambda_{\min}(X^{T}X)\|Z_{\perp}\|_{F}^{2}\label{eq:Wsol2}
\end{align}
We now plug our choice of $y$ and $W_{i,j}$ specified in (\ref{eq:ysol})
and (\ref{eq:Wsol}) into \lemref{tradeoff}, along with a good choice
of coefficients $\gamma_{1}$ and $\gamma_{2}$.
\begin{lemma}
\label{lem:costheta}We have $\cos\theta(t)\ge\alpha\cdot(t/\beta)+\sqrt{1-\alpha^{2}}\sqrt{1-t^{2}/\beta^{2}}$
for $0\le t\le\alpha\beta$ where $Z_{\perp}=(I-XX^{\dagger})Z$ and
\begin{align*}
\alpha & =\frac{\|Z_{\perp}Z_{\perp}^{T}\|_{F}}{\|XX^{T}-ZZ^{T}\|_{F}}, & \beta & =\frac{\lambda_{\min}(X^{T}X)\cdot\tr(Z_{\perp}Z_{\perp}^{T})}{\|XX^{T}-ZZ^{T}\|_{F}\|Z_{\perp}Z_{\perp}^{T}\|_{F}}.
\end{align*}
\end{lemma}
\begin{proof}
Write $\Pi=I-ZZ^{\dagger}$ and $\e=\vect(XX^{T}-ZZ^{T})$, and decompose
$\e_{1}=(\J_{X}\J_{X}^{\dagger})\e$ and $\e_{2}=(I-\J_{X}\J_{X}^{\dagger})\e$.
We verify that the following heuristic choice
\begin{gather*}
y=\underbrace{\frac{\sqrt{1-(t/\beta)^{2}}}{\|\e\|\|\e_{1}\|}}_{\gamma_{1}}\;\J_{X}^{\dagger}\e,\quad W=\underbrace{\frac{t/\beta}{\|\e\|\|\e_{2}\|}}_{\gamma_{2}}\;\sum_{i=1}^{r}\vect(Z_{\perp}e_{i}v_{r}^{T})\vect(Z_{\perp}e_{i}v_{r}^{T})^{T},
\end{gather*}
is feasible for (\ref{eq:cosprob}). Indeed, $W\succeq0$ holds, and
we have from (\ref{eq:e2}), (\ref{eq:Wsol1}), and (\ref{eq:Wsol2}):
\begin{gather*}
\f=\J_{X}y-\sum_{i=1}^{r}\vect(W_{i,i})=\gamma_{1}\cdot\J_{X}\J_{X}^{\dagger}\e-\gamma_{2}\vect(Z_{\perp}Z_{\perp}^{T})=\gamma_{1}\cdot\e_{1}+\gamma_{2}\e_{2},\\
\|\e\|^{2}\|\f\|^{2}=\gamma_{1}^{2}\cdot\|\e\|^{2}\|\e_{1}\|^{2}+\gamma_{2}^{2}\cdot\|\e\|^{2}\|\e_{2}\|^{2}=(1-(t/\beta)^{2})+(t/\beta)^{2}=1,\\
\langle\J_{X}^{T}\J_{X},W\rangle=2\left(\frac{t/\beta}{\|\e\|\|\e_{2}\|}\right)\cdot\lambda_{\min}(X^{T}X)\|Z_{\perp}\|_{F}^{2}=\frac{2t}{\beta}\cdot\frac{\sigma_{r}^{2}\cdot\tr(Z_{\perp}Z_{\perp}^{T})}{\|\e\|\|Z_{\perp}Z_{\perp}^{T}\|_{F}}=2t.
\end{gather*}
To verify that $\Pi\,\mat(\f)\,\Pi\succeq0$, note that $t\le\alpha\beta$
implies $t/\beta\le\alpha\le1$, so it follows from $\alpha=\|\e_{2}\|/\|\e\|$
and $\sqrt{1-\alpha^{2}}=\|\e_{1}\|/\|\e\|$ that
\[
\gamma_{1}=\frac{\sqrt{1-(t/\beta)^{2}}}{\|\e_{1}\|\|\e\|}\ge\frac{\sqrt{1-\alpha^{2}}}{\|\e_{1}\|\|\e\|}=\frac{1}{\|\e\|^{2}}=\frac{\alpha}{\|\e_{2}\|\|\e\|}\ge\frac{t/\beta}{\|\e_{2}\|\|\e\|}=\gamma_{2}.
\]
In turn, $\gamma_{1}\ge\gamma_{2}$ implies
\begin{align*}
\f & =\gamma_{1}\cdot\e_{1}+\gamma_{2}\cdot\e_{2}=\gamma_{1}\cdot\e-(\gamma_{1}-\gamma_{2})\e_{2},\\
\mat(\f) & =\gamma_{1}\cdot(XX^{T}-ZZ^{T})+(\gamma_{1}-\gamma_{2})\cdot Z_{\perp}Z_{\perp}^{T},\\
\Pi\,\mat(\f)\,\Pi & =\gamma_{1}\cdot\Pi(XX^{T})\Pi+(\gamma_{1}-\gamma_{2})\cdot\Pi(Z_{\perp}Z_{\perp}^{T})\Pi\succeq0.
\end{align*}
Having verified $y$ and $W$ as feasible for (\ref{eq:cosprob}),
we use its objective value to heuristically lower-bound the maximum
\[
\cos\theta(t)\ge\e^{T}\f=\gamma_{1}\cdot\|\e_{1}\|^{2}+\gamma_{2}\cdot\|\e_{2}\|^{2}=\sqrt{1-(t/\beta)^{2}}\sqrt{1-\alpha^{2}}+(t/\beta)\cdot\alpha.
\]
 \qed
\end{proof}

Finally, we evoke the following lemma, which was previously used in
the proof of \citet[Theorem~12]{zhang2019sharp} to fully solve the
rank-1 case in closed-form.
\begin{lemma}
\label{lem:optim}Given $0\le\alpha<1$ and $\beta\ge0$, let $\gamma_{\alpha,\beta}:\R_{+}\to[0,1]$
satisfy $\gamma_{\alpha,\beta}(t)=\alpha\cdot(t/\beta)+\sqrt{1-\alpha^{2}}\sqrt{1-t^{2}/\beta^{2}}$
for $t\le\alpha\beta$. Then, we have
\[
\max_{0\le t\le\alpha\beta}\frac{1+\gamma_{\alpha,\beta}(t)}{2t+1-\gamma_{\alpha,\beta}(t)}=\begin{cases}
{\displaystyle \frac{1+\sqrt{1-\alpha^{2}}}{1-\sqrt{1-\alpha^{2}}}} & \text{if }\beta\ge{\displaystyle \frac{\alpha}{1+\sqrt{1-\alpha^{2}}}},\\
{\displaystyle \frac{1-\alpha\beta}{\beta(\alpha-\beta)}} & \text{if }\beta\le{\displaystyle \frac{\alpha}{1+\sqrt{1-\alpha^{2}}}}.
\end{cases}
\]
\end{lemma}
The proof of \lemref{abdef} immediately follows by substituting \lemref{costheta}
into \lemref{tradeoff} and then solving the resulting optimization
over $t$ using \lemref{optim}.

\section{\label{sec:counterexamples}Counterexample that attains the heuristic
upper-bound (\lemref{heuristic})}

Finally, we give a constructive proof for the heuristic upper-bound
in \lemref{heuristic}. First, recall from \lemref{sdp_formu} that
our SDP upper-bound is defined as the following
\begin{equation}
\kappa_{\ub}(X,Z)=\min_{\kappa,\H}\left\{ \kappa:I\preceq\H\preceq\kappa I,\:\J_{X}^{T}\H\e=0,\;-2I_{r}\otimes\mat(\H\e)\preceq\J_{X}^{T}\H\J_{X}\right\} \label{eq:kappaub_prob}
\end{equation}
where $\e=\vect(XX^{T}-ZZ^{T})$ and its Jacobian $\J_{X}$ is defined
to satisfy $\J_{X}\vect(V)=\vect(XV^{T}+VX^{T})$ for all $V$. Below,
we derive a set of sufficient conditions for a given $(\kappa,\H)$
to be feasible for (\ref{eq:kappaub_prob}). 
\begin{lemma}
\label{lem:counter}Let $[Q_{1},Q_{2}]$ have orthonormal columns
with $Q_{1}\in\R^{n\times r}$ and $Q_{2}\in\R^{n\times r^{\star}}$.
If $\kappa=1+2\sqrt{r/r^{\star}}$, and $\H$ is any $n^{2}\times n^{2}$
matrix that satisfies $I\preceq\H\preceq\kappa I$ and the following
$rr^{\star}+2$ eigenvalue equations
\begin{gather*}
\H\vect\left(Q_{1}Q_{1}^{T}+\sqrt{r/r^{\star}}Q_{2}Q_{2}\right)=\kappa\cdot\vect\left(Q_{1}Q_{1}^{T}+\sqrt{r/r^{\star}}Q_{2}Q_{2}\right),\\
\H\vect\left(Q_{1}Q_{1}^{T}-\sqrt{r/r^{\star}}Q_{2}Q_{2}\right)=1\cdot\vect\left(Q_{1}Q_{1}^{T}-\sqrt{r/r^{\star}}Q_{2}Q_{2}\right),\\
\H\vect\left(Q_{1}VQ_{2}^{T}+Q_{2}V^{T}Q_{1}^{T}\right)=\kappa\cdot\vect\left(Q_{1}VQ_{2}^{T}+Q_{2}V^{T}Q_{1}^{T}\right)\qquad\text{for all }V\in\R^{r\times r^{\star}},
\end{gather*}
then $(\kappa,\H)$ is a feasible point for (\ref{eq:kappaub_prob})
with $X=Q_{1}$ and $Z=\sqrt{1+\sqrt{r/r^{\star}}}Q_{2}$. Moveover,
$\e^{T}\H\e=(1+2\sqrt{r/r^{\star}})\cdot(1+\sqrt{r/r^{\star}})\cdot r^{\star}$.
\end{lemma}
\begin{proof}
Write $\rho\equiv\sqrt{r/r^{\star}}$, and take $[Q_{1},Q_{2}]=I_{r+r^{\star}}$
without loss of generality. Decomposing $\e=\vect(XX^{T}-ZZ^{T})$
into eigenvectors and applying $\H$ yields
\begin{gather*}
\mat(\e)=\begin{bmatrix}I_{r} & 0\\
0 & -(1+\rho)I_{r^{\star}}
\end{bmatrix}=-\frac{1}{2\rho}\begin{bmatrix}I_{r} & 0\\
0 & \rho I_{r^{\star}}
\end{bmatrix}+\frac{1+2\rho}{2\rho}\begin{bmatrix}I_{r} & 0\\
0 & -\rho I_{r^{\star}}
\end{bmatrix},\\
\mat(\H\e)=-\frac{1+2\rho}{2\rho}\begin{bmatrix}I_{r} & 0\\
0 & \rho I_{r^{\star}}
\end{bmatrix}+\frac{1+2\rho}{2\rho}\begin{bmatrix}I_{r} & 0\\
0 & -\rho I_{r^{\star}}
\end{bmatrix}=\begin{bmatrix}0_{r} & 0\\
0 & -\kappa I_{r^{\star}}
\end{bmatrix}.
\end{gather*}
Clearly, $\inner{\mat(\e)}{\mat(\H\e)}=\kappa(1+\rho)r^{\star}.$
For $v=\vect(V)=\vect([V_{1};V_{2}])$ with $V_{1}\in\R^{r\times r}$
and $V_{2}\in\R^{r^{\star}\times r}$, observe that 
\begin{gather*}
v^{T}\J_{X}\H\e=\left\langle \begin{bmatrix}V_{1}+V_{1}^{T} & V_{2}^{T}\\
V_{2} & 0_{r^{\star}}
\end{bmatrix},\begin{bmatrix}0_{r} & 0\\
0 & -\kappa I_{r^{\star}}
\end{bmatrix}\right\rangle =0,\\
v^{T}[I_{r}\otimes\mat(\H\e)]v=\langle\mat(\H\e),VV^{T}\rangle=\left\langle \begin{bmatrix}0_{r} & 0\\
0 & -\kappa I_{r^{\star}}
\end{bmatrix},\begin{bmatrix}V_{1}V_{1}^{T} & V_{1}V_{2}^{T}\\
V_{2}V_{1}^{T} & V_{2}V_{2}^{T}
\end{bmatrix}\right\rangle =-\kappa\|V_{2}\|^{2},\\
v^{T}\J_{X}^{T}\H\J_{X}v=\left\Vert \H^{1/2}\vect\left(\begin{bmatrix}V_{1}+V_{1}^{T} & V_{2}^{T}\\
V_{2} & 0_{r^{\star}}
\end{bmatrix}\right)\right\Vert ^{2}\ge\left\Vert \H^{1/2}\vect\left(\begin{bmatrix}0_{r} & V_{2}^{T}\\
V_{2} & 0_{r^{\star}}
\end{bmatrix}\right)\right\Vert ^{2}=2\kappa\|V_{2}\|^{2}.
\end{gather*}
Therefore, $v^{T}\J_{X}\H\e=0$ and $2v^{T}[I_{r}\otimes\mat(\H\e)]v\ge-v^{T}\J_{X}^{T}\H\J_{X}v$
hold for all $v$. \qed
\end{proof}

We now state an explicit counterexample that satisfies the conditions
in \lemref{counter}. Without loss of generality, take $u_{1},\dots,u_{n}$
to be the standard basis for $\R^{n}$. The counterexample is constructed
by taking the rescaled standard basis $A^{(i,j)}=\sqrt{\kappa}u_{i}u_{j}^{T}$
for $\R^{n\times n}$, and then replacing the elements $\{A^{(i,i)}\}_{i=1}^{n}$
with a new orthogonal basis for $\diag(\R^{n})$ whose first two elements
$A^{(1,1)}$ and $A^{(2,2)}$ are exactly the first two eigenvectors
in Lemma 6.1. The ordering $\pi$ for the Gram--Schmidt orthogonalization
below is chosen specifically to ensure that $\mathrm{span}\{A^{(i,i)}\}_{i=1}^{n}=\diag(\R^{n})$. 
\begin{example}[$\kappa=1+2\sqrt{r/r^{\star}}$ counterexample]
\label{exa:overparam}Given $r,r^{\star},n$ satisfying $1\le r^{\star}\le r<n$,
let $u_{1},u_{2},\ldots,u_{n}$ be an orthonormal basis for $\R^{n}$.
Define the quadratic function
\[
\phi_{0}(M)=\frac{1}{2}\sum_{i=1}^{n}\sum_{j=1}^{n}\left|\langle A^{(i,j)},M-M^{\star}\rangle\right|^{2}\quad\text{ where }M^{\star}=\sum_{k=r+1}^{r+r^{\star}}u_{k}u_{k}^{T},
\]
and $A^{(i,j)}=\sqrt{\kappa}\cdot u_{i}u_{j}^{T}$ for all $1\le i,j\le n$
except the following
\begin{gather*}
A^{(1,1)}=\sqrt{\frac{\kappa}{2r}}\sum_{j=1}^{r}u_{j}u_{j}^{T}+\sqrt{\frac{\kappa}{2r^{\star}}}\sum_{k=r+1}^{r+r^{\star}}u_{k}u_{k}^{T},\\
A^{(2,2)}=\frac{1}{\sqrt{2r}}\sum_{j=1}^{r}u_{j}u_{j}^{T}-\frac{1}{\sqrt{2r^{\star}}}\sum_{k=r+1}^{r+r^{\star}}u_{k}u_{k}^{T},\\
A^{(i,i)}=u_{\pi(i)}u_{\pi(i)}^{T}-\sum_{j=1}^{i-1}\frac{A^{(j,j)}}{\|A^{(j,j)}\|_{F}^{2}}\inner{A^{(j,j)}}{u_{\pi(i)}u_{\pi(i)}^{T}}\text{ for }i\in\{3,\dots,r+r^{\star}\},
\end{gather*}
with the ordering $\pi=(1,r+1,2,3,\dots,r,r+2,r+3,\dots,r+r^{\star})$.
Note that $\{A^{(i,j)}\}$ is an orthogonal basis for $\R^{n\times n}$
with $1\le\|A^{(i,j)}\|_{F}^{2}\le\kappa$, so the quadratic function
$\phi_{0}$ is twice-differentiable, $1$-strongly convex and $\kappa$-smooth,
\[
\|E\|_{F}^{2}\le\inner{\nabla^{2}\phi_{0}(M)[E]}E\le\kappa\|E\|_{F}^{2}\qquad\text{for all }M\in\R^{n\times n}.
\]
Its minimizer $M^{\star}=\arg\min_{M}\phi_{0}(M)$ satisfies $r^{\star}=\rank(M^{\star}),$
but its corresponding $f_{0}(U)\eqdef\phi_{0}(UU^{T})$ admits the
$n\times r$ matrix $X=\frac{1}{\sqrt{1+\sqrt{r/r^{\star}}}}[u_{1},\dots,u_{r}]$
as a second-order critical point 
\begin{align*}
f_{0}(X)-\min_{U}f_{0}(U) & =\left(1+\frac{1}{1+\sqrt{r^{\star}/r}}\right)\frac{\|M^{\star}\|_{F}^{2}}{2}, & \nabla f_{0}(X) & =0, & \nabla^{2}f_{0}(X) & \succeq0.
\end{align*}
 \qed
\end{example}

\section{\label{sec:strictsaddle}Proof of the strict saddle property (\propref{strictsaddle})}

\global\long\def\Set{\mathcal{S}}%
We now prove that, in our setting of an $L$-smooth and $\mu$-strongly
convex $\phi$, the absence of spurious local minima in $f(U)=\phi(UU^{T})$
implies the \emph{strict saddle property}~\citep{ge2015escaping,jin2017escape,ge2017nospurious}.
Let $\mathcal{S}(\epsilon,\rho)$ denote the set of spurious $\epsilon$
second-order critical points for its corresponding Burer--Monteiro
function $f(U)\eqdef\phi(UU^{T})$ that lie more than a distance of
$\rho$ from the minimizer:
\[
\Set(\epsilon,\rho)=\{X\in\R^{n\times r}:\|\nabla f(X)\|_{F}\le\epsilon,\quad\nabla^{2}f(X)\succeq-\epsilon I,\quad\|XX^{T}-M^{\star}\|_{F}\ge\rho\}.
\]
Our crucial observation is that $\Set(\epsilon,\rho)$ is a \emph{compact}
set: (i) it is by definition closed; (ii) it is bounded because any
$X\in\Set(\epsilon,\rho)$ with $\|X\|_{F}\to\infty$ would imply
$\|\nabla f(X)\|_{F}\to\infty$ via \Lemref{grad_coercive} below. 
\begin{lemma}
\label{lem:grad_coercive}Let $\phi$ be $L$-smooth and $\mu$-strongly
convex, and let $M^{\star}=\arg\min_{M\succeq0}\phi(M)$ be attained.
Then, the gradient norm $\|\nabla f(X)\|_{F}$ of the Burer--Monteiro
function $f(X)\eqdef\phi(XX^{T})$ is a coercive function, meaning
that $\|\nabla f(X)\|_{F}\to\infty$ as $\|X\|_{F}\to\infty$.
\end{lemma}
\begin{proof}
It follows from the $L$-smoothness and $\mu$-strong convexity of
$\phi$ that its Hessian approximately preserves inner products (see
e.g.~\citet[Proposition~2.1]{li2019non})
\[
\frac{2}{\mu+L}\inner{\nabla^{2}\phi(M)[E]}F\ge\inner EF-\frac{L-\mu}{L+\mu}\|E\|_{F}\|F\|_{F}.
\]
Moreover, the fact that $M^{\star}=\arg\min_{M\succeq0}\phi(M)$ implies
that $\nabla\phi(M^{\star})\succeq0$. It follows from the fundamental
theorem of calculus that
\begin{align*}
\inner{\nabla f(X)}X & =2\inner{\nabla\phi(XX^{T})}{XX^{T}}\\
 & =2\inner{\int_{0}^{1}\nabla^{2}\phi(M_{t})[XX^{T}-M^{\star}]\,\d t}{XX^{T}}+2\inner{\nabla\phi(M^{\star})}{XX^{T}}\\
 & \ge(L+\mu)\inner{XX^{T}-M^{\star}}{XX^{T}}-(L-\mu)\|XX^{T}-M^{\star}\|_{F}\|XX^{T}\|_{F}\\
 & \ge2\mu\|XX^{T}\|_{F}^{2}-2L\|XX^{T}\|_{F}\|M^{\star}\|_{F}\ge2\|X\|_{F}^{2}\left(\frac{\mu}{r}\|X\|_{F}^{2}-L\|M^{\star}\|_{F}\right)
\end{align*}
where we used $\sqrt{r}\|XX^{T}\|_{F}\ge\tr XX^{T}$ for $X\in\R^{n\times r}$
and $n\ge r$. Therefore, once $\|X\|_{F}^{2}>\frac{Lr}{\mu}\|M^{\star}\|_{F}$
holds, we have $\|\nabla f(X)\|_{F}\|X\|_{F}\ge\inner{\nabla f(X)}X\ge C\|X\|_{F}^{2}$
for some $C>0$. \qed
\end{proof}

The strict saddle property then follows immediately from compactness.\footnote{This proof was first suggested by a reviewer of the paper~\citep{zhang2021preconditioned}
at \url{https://openreview.net/forum?id=5-Of1DTlq&noteId=FWahQu0GxmD}} 
\begin{lemma}
\label{lem:saddle}Let $\phi:\S^{n}\to\R$ be twice-differentiable,
$L$-smooth and $\mu$-strongly convex, and let the minimizer $M^{\star}=\arg\min_{M\succeq0}\phi(M)$
exist. If the function $f(U)\eqdef\phi(UU^{T})$ satisfies the following
\[
\nabla f(X)=0,\quad\nabla^{2}f(X)\succeq0\quad\iff\quad XX^{T}=M^{\star}
\]
then it also satisfies the strict saddle property 
\[
\|\nabla f(\tilde{X})\|_{F}\le\epsilon(\rho),\quad\nabla^{2}f(\tilde{X})\succeq-\epsilon(\rho)I\quad\implies\quad\|\tilde{X}\tilde{X}^{T}-M^{\star}\|_{F}<\rho
\]
where the function $\epsilon$ is nondecreasing and satisfies $0<\epsilon(\rho)\le1$
for all $\rho>0$.
\end{lemma}
\begin{proof}
For an arbitrary $\rho>0$, define the sequence of compact subsets
$\Set_{1}\supseteq\Set_{2}\supseteq\Set_{3}\supseteq\cdots$ where
$\Set_{k}=\Set(\frac{1}{k},\rho)$. Given that $f$ has no spurious
local minima, the intersection of all $\{\Set_{k}\}_{k=1}^{\infty}$
must be empty. It follows from the finite open cover property of the
compact set $\Set_{1}$ that there exists some finite $N\ge1$ such
that $\Set_{N}=\emptyset$. Therefore, any $X$ that satisfies $\|\nabla f(X)\|_{F}\le\frac{1}{N}$
and $\nabla^{2}f(X)\succeq-\frac{1}{N}I$ must also satisfy $\|XX^{T}-M^{\star}\|_{F}<\rho$.
Finally, we repeat this for all $\rho>0$ and define $\epsilon(\rho)=\frac{1}{N}$.
We observe that the function $\epsilon(\cdot)$ satisfies $\epsilon(\rho)\in(0,1]$
because $N\ge1$ is finite, and that it is nondecreasing because $\Set_{k+1}\subseteq\Set_{k}$
for all $k$. \qed
\end{proof}

We conclude the proof by using the convexity and $L$-smoothness of
$\phi$ to show that an $\epsilon$-second-order point $X$ that is
also $\rho$-close is guaranteed to be $\delta$-globally suboptimal.
\begin{proof}
[\propref{strictsaddle}]Write $M^{\star}=ZZ^{T}$, and let $X\in\R^{n\times r}$
be an $\epsilon$-second-order point satisfying $\|\nabla f(X)\|\le\epsilon$
and $\nabla^{2}f(X)\succeq-\epsilon I$. If $f$ has no spurious local
minima, then \lemref{saddle} says that $X$ is also $\rho(\epsilon)$-close,
as in $\|XX^{T}-ZZ^{T}\|_{F}<\rho(\epsilon)$ where $\rho$ is a nondecreasing
function such that $\rho(\epsilon)\to0$ as $\epsilon\to0$. If overparameterized
$r>\rank(Z)$, we use the rank deficiency of $X$ to certify global
optimality \citep[Proposition~10]{zhang2022preconditioned}
\begin{align*}
f(X)-f(Z) & \le\frac{1}{2}\|X\|_{F}\|\nabla f(X)\|_{F}-\frac{1}{2}\|Z\|_{F}^{2}\lambda_{\min}(\nabla^{2}f(X))+2L\|Z\|_{F}^{2}\lambda_{\min}(X^{T}X)\\
 & \le\left[\frac{1}{2}(\|Z\|_{F}+\rho)+\frac{1}{2}\|Z\|_{F}^{2}\right]\epsilon+2L\|Z\|_{F}^{2}\rho(\epsilon)\eqdef\varphi(\epsilon).
\end{align*}
Note that the near rank deficiency of $X$ is due to Weyl's inequality
$\lambda_{r}(XX^{T})\le\lambda_{r}(ZZ^{T})+\|XX^{T}-ZZ^{T}\|_{F}\le0+\rho(\epsilon)$.
The function $\varphi(\epsilon)$ is nondecreasing and satisfies $\varphi(\epsilon)\to0$
as $\epsilon\to0$. We define $\epsilon(\delta)=\varphi^{-1}(\delta)$
to ensure that an $\epsilon(\delta)$-second-order point $X$ will
satisfy $f(X)-f(Z)\le\delta$.

In the exactly parameterized regime $r=\rank(Z)$, we use the displacement
vector $\Delta=X-ZR$ where $R=\arg\min_{RR^{T}=R^{T}R=I}\|X-ZR\|_{F}$;
see also~\citep{tu2016low,ge2017nospurious,zhu2018global,li2019non,chi2019nonconvex}.
Within a small neighborhood $\rho(\epsilon)\le\frac{1}{2}\lambda_{\min}(Z^{T}Z)$,
the norm of the displacement vector scales with the error norm \citep[Lemma~5.4]{tu2016low}:
\[
\|\Delta\|_{F}^{2}\le\frac{\|XX^{T}-ZZ^{T}\|_{F}^{2}}{2(\sqrt{2}-1)\lambda_{\min}(X^{T}X)}\le\frac{\rho(\epsilon)^{2}}{\frac{1}{2}[\lambda_{\min}(Z^{T}Z)-\rho(\epsilon)]}\le\frac{4\rho(\epsilon)^{2}}{\lambda_{\min}(Z^{T}Z)}.
\]
The error vector decomposes as $XX^{T}-ZZ^{T}=X\Delta+\Delta X^{T}-\Delta\Delta^{T}$,
and this yields via the convexity of $\phi$:
\begin{align*}
f(X)-f(Z) & =\phi(XX^{T})-\phi(ZZ^{T})\le\inner{\nabla\phi(XX^{T})}{XX^{T}-ZZ^{T}}\\
 & =\inner{\nabla\phi(XX^{T})}{X\Delta+\Delta X^{T}}-\inner{\nabla\phi(XX^{T})}{\Delta\Delta^{T}}\\
 & \le\inner{\nabla f(X)}{\Delta}-\frac{1}{2}\left[\inner{\nabla^{2}f(X)[\Delta]}{\Delta}-L\|X\Delta^{T}+\Delta X^{T}\|_{F}^{2}\right]\\
 & \le\epsilon\|\Delta\|_{F}+\frac{1}{2}\epsilon\|\Delta\|_{F}^{2}+L\lambda_{\max}(X^{T}X)\|\Delta\|_{F}^{2}\\
 & \le\left[\frac{2\epsilon}{\sqrt{\lambda_{\min}(Z^{T}Z)}}\right]\rho(\epsilon)+\left[\frac{2\epsilon+4L[\lambda_{\max}(Z^{T}Z)+\rho(\epsilon)]}{\lambda_{\min}(Z^{T}Z)}\right]\rho(\epsilon)^{2}\eqdef\varphi(\epsilon).
\end{align*}
The function $\varphi(\epsilon)$ is nondecreasing and satisfies $\varphi(\epsilon)\to0$
as $\epsilon\to0$. We define $\epsilon(\delta)=\min\{\varphi^{-1}(\delta),\rho^{-1}(\frac{1}{2}\lambda_{\min}(Z^{T}Z))\}$
to ensure that an $\epsilon(\delta)$-second-order point $X$ will
satisfy $f(X)-f(Z)\le\delta$. \qed
\end{proof}

\section*{Acknowledgments}

I thank Salar Fattahi and Cedric Josz for introducing me to overparameterization,
and also Sabrina Zielinski, Gavin Zhang, Nicolas Boumal, Iosif Pinelis,
Nati Srebro, and the three anonymous reviewers for helpful discussions
and feedback. I am grateful to Salar Fattahi for connecting Theorem~\ref{thm:EckartYoung}
with the Eckart--Young Theorem; this turned out to be a critical
insight for the paper. The proof of \lemref{saddle} is due to an
anonymous reviewer of the paper~\citep{zhang2021preconditioned}.
The proof of \thmref{EckartYoung} has been substantially simplified
by Reviewer 2 via a trace inequality of Ruhe. I thank Reviewer 3 for
motivating me to re-examine \exaref{overparam}, and Nati Srebro for
prompting me to close the gap between necessity and sufficiency in
\thmref{main}.

\subsubsection*{Funding}

Financial support for this work was provided in part by the NSF CAREER
Award ECCS-2047462 and in part by C3.ai Inc. and the Microsoft Corporation
via the C3.ai Digital Transformation Institute.

\subsubsection*{Conflicts of interest/Competing interests}

The author has no relevant financial or non-financial interests to
disclose.

\appendix

\section{\label{app:tradeoff}Derivation of the dual for the lower-bound problem
(Proof of \lemref{tradeoff})}

We will derive the Lagrangian dual for the following
\begin{align*}
\kappa_{\lb}(X,Z) & =\min_{\kappa,\H,s}\left\{ \kappa:\begin{array}{c}
I\preceq\H\preceq\kappa I,\qquad\J_{Z}^{T}s=0,\qquad s\in\vect(\S_{+}^{n}),\\
\J_{X}^{T}(\H\e+s)=0,\quad2I_{r}\otimes\mat(\H\e+s)+\kappa\J_{X}^{T}\J_{X}\succeq0.
\end{array}\right\} 
\end{align*}
Recall that $\e=\vect(XX^{T}-ZZ^{T})$ and $\J_{X}\vect(V)=\vect(XV^{T}+VX^{T})$
with respect to fixed $X,Z\in\R^{n\times r}$. Observe that
\[
\J_{Z}^{T}\vect(S)=0\iff(S+S^{T})Z=0\iff S=Q_{\perp}S_{\perp}Q_{\perp}^{T}\iff\vect(S)=\Q_{\perp}\vect(S_{\perp})
\]
where the orthogonal complement $Q_{\perp}$ of $Z$ is such that
$Q_{\perp}Q_{\perp}^{T}=I_{n}-ZZ^{\dagger}$, and $\Q_{\perp}=Q_{\perp}\otimes Q_{\perp}$.
Substituting and taking the Lagrangian dual yields
\begin{align}
\kappa_{\lb}(X,Z) & =\min_{\kappa,\H,s_{\perp}}\left\{ \kappa:\begin{array}{c}
I\preceq\H\preceq\kappa I,\quad s_{\perp}\in\vect(\S_{+}^{n-r^{\star}}),\quad\J_{X}^{T}(\H\e+\Q_{\perp}s_{\perp})=0,\\
2I_{r}\otimes\mat(\H\e+\Q_{\perp}s_{\perp})+\kappa\J_{X}^{T}\J_{X}\succeq0
\end{array}\right\} \nonumber \\
 & =\max_{y,W_{i,j},U,V}\left\{ \tr(U):\begin{array}{c}
\f\e^{T}+\e\f^{T}=U-V,\quad\tr(V)+\inner{\J_{X}^{T}\J_{X}}W=1,\\
\f\equiv\J_{X}y-\sum_{i=1}^{r}\vect(W_{i,i}),\quad\Q_{\perp}^{T}\f\in\vect(\S_{+}^{n-r^{\star}}),\\
U,V\succeq0,\quad W=[W_{i,j}]_{i,j=1}^{r}\succeq0
\end{array}\right\} \label{eq:ub2-1}
\end{align}
in which $y\in\R^{nr}$ and $W_{i,j}\in\R^{n\times n}$ for $i,j\in\{1,2,\dots,r\}$.
Strong duality holds here because $y=0,$ $W=\epsilon\cdot I_{nr},$
$U=\frac{1}{n}I_{n},$ and $V=\frac{1}{n}I_{n}+\epsilon\cdot r\cdot[\e\vect(I_{n})^{T}+\e\vect(I_{n})^{T}]$
is strictly feasible for a sufficiently small $\epsilon>0$. 

The dual problem (\ref{eq:ub2-1}) has a closed-form solution over
$U$ and $V$ \citep[Lemma~13]{zhang2019sharp}
\begin{align*}
\frac{\tr([M]_{-})+\alpha}{\tr([M]_{+})}= & \min_{t,U,V}\left\{ \tr(V)+\alpha\cdot t:\begin{array}{c}
tM=U-V,\quad t\ge0,\\
\tr(U)=1,\quad U,V\succeq0
\end{array}\right\} ,
\end{align*}
and $\f\e^{T}+\e\f^{T}$ has exactly two nonzero eigenvalues $\e^{T}\f\pm\|\e\|\|\f\|$
\citep[Lemma~14]{zhang2019sharp}. Substituting this solution yields
exactly 
\begin{align*}
\kappa_{\lb}(X,Z) & =\max_{y,W_{i,j}}\left\{ \frac{\|\e\|\|\f\|+\e^{T}\f}{\|\e\|\|\f\|-\e^{T}\f+\inner{\J_{X}^{T}\J_{X}}W}:\begin{array}{c}
\Q_{\perp}\Q_{\perp}^{T}\f\in\vect(\S_{+}^{n}),\\
\f=\J_{X}y-\sum_{i=1}^{r}\vect(W_{i,i}),\\
W=[W_{i,j}]_{i,j=1}^{r}\succeq0
\end{array}\right\} \\
 & =\max_{t\ge0}\frac{1+\cos\theta(t)}{2t+1-\cos\theta(t)}
\end{align*}
where $\cos\theta(t)$ is itself is defined as the following maximization
\[
\cos\theta(t)=\max_{y,W_{i,j}}\left\{ \e^{T}\f:\begin{array}{c}
\f\equiv\J_{X}y-\sum_{i=1}^{r}\vect(W_{i,i}),\quad\|\e\|\|\f\|=1,\\
(I-ZZ^{\dagger})\,\mat(\f)\,(I-ZZ^{\dagger})\succeq0,\\
\langle\J_{X}^{T}\J_{X},W\rangle=2t,\quad W=[W_{i,j}]_{i,j=1}^{r}\succeq0.
\end{array}\right\} 
\]

\bibliographystyle{abbrvnat}
\bibliography{proof_half}

\end{document}